   \def\obs#1{{\bf (*** #1 ***)} }

 \def\obs#1{}     

\documentclass[twoside,letterpaper,draft,11pt]{amsart}

\usepackage{amsmath}               
\usepackage{amsthm}                
\usepackage{latexsym}

\usepackage{xspace}
\usepackage{amscd}

\title[Crossed products by twisted partial actions and graded algebras]{Crossed products by
twisted partial actions and graded algebras}

\thanks{{This work was partially supported by CNPq of Brazil and
Secretar\'{\i}a de Estado de Universidades e Investigaci\' on del MEC, Espa\~ na.
}}

\newtheorem{teo}{Theorem}[section]
\newtheorem{defi}[teo]{Definition}
\newtheorem{lema}[teo]{Lemma}
\newtheorem{cor}[teo]{Corollary}
\newtheorem{prop}[teo]{Proposition}

\theoremstyle{remark}\newtheorem{remark}[teo]{Remark}
\theoremstyle{remark}\newtheorem{ex}[teo]{Example}

\font\msbm = msbm10
\newcommand{\stimes}{{\hbox {\msbm o}}}

\newcommand{\Cc}{{\mathbb C}}

\newcommand{\Rr}{{\mathbb R}}

\newcommand{\Nn}{{\mathbb N}}

\newcommand{\0}{\theta}

\newcommand{\D}{{\mathcal D}}

\newcommand{\de}{\delta}

\newcommand{\af}{\alpha}

\newcommand{\z}{\zeta}

\newcommand{\p}{{\bf Proof. }}

\newcommand{\fim}{\hfill\mbox{$\Box$}}

\newcommand{\C}{{\mathcal C}}
\newcommand{\A}{{\mathcal A}}
\newcommand{\R}{{\mathcal R}}
\newcommand{\B}{{\mathcal B}}
\newcommand{\I}{{\mathcal I}}
\newcommand{\J}{{\mathcal J}}
\newcommand{\M}{{\mathcal M}}
\newcommand{\K}{{\mathcal K}}

\newcommand{\la}{{\longrightarrow}}

\newcommand{\m}{{}^{-1}}
\newcommand{\mt}{\mapsto}
\def\f{\varphi}
\def\e{\varepsilon}

\def\ndv{\ {\mid \kern -0.7 em {\scriptstyle \not}} \ \ }

\def\nd{\ {\mid \kern -0.4 em {\scriptstyle \not}} \ \ }

\newcommand{\Mod}{\mbox{\rm -mod}}

\newcommand{\Hom}{\mathrm{Hom}}
\newcommand{\N}{{\mathbb N}}
\newcommand{\E}{{\mathcal E}}

\begin{document}

\date{\today}

\maketitle

\begin{center}
{\large M. Dokuchaev}\\ {\footnotesize Instituto de
Matem\'atica e Estat\'\i stica\\
Universidade de S\~ao Paulo, \\
05508-090 S\~ao Paulo, SP, Brasil
\\ E-mail: {\it dokucha@ime.usp.br}}
\end{center}

\begin{center}
{\large R.\ Exel}\\ {\footnotesize Departamento de Matem\'{a}tica,\\
 Universidade Federal de Santa Catarina,\\
 88040-900 Florian\'opolis, SC,  Brasil\\
E-mail:\it{exel@mtm.ufsc.br}}
\end{center}

\begin{center}
{\large J.\ J.\ Sim\'{o}n}\\  {\footnotesize Departamento de Matem\'{a}ticas,\\
 Universidad de Murcia,\\
 30071 Murcia, Espa\~{n}a\\
E-mail:\it{jsimon@um.es}}
\end{center}

\begin{abstract} For a twisted partial action $\Theta $ of a group $G$ on
an (associative  non-necessarily unital) algebra  $ \A$ over a
commutative unital ring $k,$ the crossed product $\A
\stimes_{\Theta} G$ is proved to be associative. Given
 a $G$-graded $k$-algebra $\B = \oplus _{g \in G} \B _g$  with
the mild restriction of homogeneous non-degeneracy, a criteria is
established for $\B $ to be  isomorphic to the crossed product $\B
_1 \stimes_{\Theta} G$ for some twisted partial action of $G$ on
$\B _1.$ The equality $B_g \B _{g\m } B _g = \B _g \; (\forall
g\in G)$ is one of the ingredients of the criteria, and if it
holds and, moreover, $\B $ has enough local units, then it is
shown that $\B $ is stably isomorphic  to a crossed product by a
twisted partial action of $G.$
\end{abstract}

\begin{section}{Introduction}\label{sec:intro}

In the Theory of Operator Algebras partial actions of groups
appeared as  a general  approach to study $C^*$-algebras generated
by partial isometries, several  relevant classes of which turned
out to be crossed products by partial actions. Crossed products
classically are in the center of the rich interplay between
dynamical systems and operator algebras, and the efforts in
generalizing them produce  structural knowledge on algebras
underlied by the new  constructions. Thus the notion of a
$C^*$-crossed product by a partial automorphism, given in
\cite{E-1},  permitted to characterize the approximately finite
$C^*$-algebras, the Bunce-Deddence  and the
Bunce-Deddence-Toeplitz algebras   as such crossed products
\cite{E-2}, \cite{E-3}. Further generalizations and steps were
made in \cite{Mc}, \cite{E1}, \cite{E0}, \cite{E2}. In \cite{ELQ}
a machinery was developed, based on the interaction between
partial actions and partial representations, enabling to study
representations and the ideal structure of partial crossed
products, and including into consideration prominent examples such
as the Toeplitz $C^*$-algebars of quasi-lattice ordered groups, as
well as the the Cuntz-Krieger algebras. The technique was also
used in \cite{EL} to define and study Cuntz-Krieger algebras for
arbitrary infinite matrices, and recently these ideas were also
applied  to Hecke algebras defined by subnormal subgroups of
length two  \cite{E3}, as well as  to
$C^*$-algebras generated by crystals and  quasi-crystals \cite{Boava}.\\

The  general notion of a (continuous)  twisted partial action of a
locally compact group on a $C^*$-algebra (a twisted partial
$C^*$-dynamical system) and the corresponding crossed products
were introduced in \cite{E0}, where the associativity of crossed
products  was proved by means of approximate identities. The
construction permitted to show in particular that any second
countable $C^*$-algebraic bundle\footnote{A $C^*$-algebraic bundle
is roughly a ``$C^*$-algebra'' graded by a locally compact group.}
with stable unit fibre is
a crossed product by a twisted partial action \cite{E0}.\\

Algebraic counterparts for some of the above mentioned notions
were introduced and studied in  \cite{DEP} and \cite{DE},
stimulating further investigations in \cite{BCFP},
\cite{CaeDeGroot}, \cite{C}, \cite{CF}, \cite{CFM}, \cite{DFP},
\cite{DP}, \cite{DdRS}, \cite{DZh}, \cite{F},  
\cite{FL} and \cite{Jonas}.
Interesting results have been obtained also in \cite{KL} and
\cite{S} (with respect to the topological part  of \cite{KL} see
also \cite{AbadieTwo}). In \cite{KL} the authors point out several
old and more recent results in various areas of mathematics
underlied by partial group actions, and, more generally, the
importance of the partial symmetry is discussed in \cite{L}, while
developing the theory of
inverse semigroups as a natural approach to  such symmetries.\\

The goal of the present article is to pursue the algebraic part of the
programm by introducing the twisted partial actions of groups on
abstract algebras over commutative rings and the corresponding crossed
products, and comparing them with general graded algebras.  The reader
will note a change in our treatment compared to previous 
ones:
partial actions on rings, as well as on other structures, were defined
so far in terms of intersections of the involved domains, and in the
present paper we prefer to do it in terms of products. Moreover, the
domains are supposed to be idempotent and pairwise commuting ideals.
Observe firstly that this does not mean a change with respect to
$C^*$-algebras: in that context the domains are closed two-sided
ideals, and a product of such two ideals coincides with their
intersection, in particular, they commute and are
idempotent. Secondly, in the majority of the algebraic situations
considered so far, the domains are unital (or more generally
$s$-unital) rings and, as a consequence, again the products coincide
with the intersections. Thirdly, the ``product definition'' serves our
treatment of graded rings, whereas it is not even clear to us what
would be the ``nicest'' definition in terms of intersections.\\

The definition of a twisted partial action incorporates a
``partial $2$-cocycle equality'' which is definitely an
associativity ingredient. The above mentioned conditions on the
domains, i. e. to be idempotent and pairwise commuting, also serve
the associativity purpose. In Section~\ref{assoc}, after giving
our main definitions, we prove the associativity of the crossed
product in Theorem~\ref{assoc}. This is done by means of a
commuting property of left
and right multipliers (\ref{L-R}), introduced in \cite{DE} to serve the non-twisted case.\\

 Our next question is to give a
criteria for a graded algebra to be a twisted partial crossed
product. It is an easy exercise that for a group $G,$ a $G$-graded
unital algebra $\B = \bigoplus_{g\in G} B _g$ is a global crossed product
$\B _1 \stimes G$ exactly when each $\B _g $ contains an element which
is
invertible in $\B ,$ i. e.
\begin{equation*}\label{global}  \forall g \in G \; \; \exists \;  u_g \in \B _g , v_g \in \B _{g\m } \mbox{ such that }
u_g v_g = v_g u_g
=1_{\B }.
\end{equation*} The corresponding criteria with twisted partial actions is more involved: using some analogue of the
$u_g$'s and $v_g$'s one needs to construct isomorphisms $\D _{g\m
} \to \D _g$ between some ideals $\D _g$ of $\B _1,$ a ``partial
twisting'' formed by certain multipliers and everything should be
combined into  a twisted partial action $\Theta ,$  so that the
given $G$-graded algebra  $\B = \bigoplus _{g \in G} \B _g $ will
be $\B _1 \stimes G.$ The easy part is to determine  the $\D
_g$'s: simply looking at the  partial crossed product, one readily
comes to  define: $\D _{g} =  \B _{g }\B _{g \m } =\{ \sum a_i b_i
: a_i \in \B _{g\m }, b_i \in \B _g\}.$ Next  a crucial point
comes: each
  $\B _g$ is a $(\D _g, \D _{g\m })$-bimodule,   and
one immediately forms the surjective Morita context $(\D _g, \D
_{g\m }, \B _g, \B _{g\m })$ in which the bimodule  maps are
determined by the product in $\B .$ The rings in the context are
non-necessarily unital and if we start with   an arbitrary
$G$-graded $\B ,$ they may be non-idempotent and even degenerate
(in the sense of Section~\ref{sec:multipliercross}). So a
restriction which guarantees that the $\D _g$' are idempotent is
needed  and 
there
is also a reasonable assumption designed to avoid
degeneracy. The latter is provided by the homogeneous
non-degeneracy condition (\ref{hnondegenerate}) (see
Section~\ref{sec:multipliercross} ) which we impose on our graded
algebra\footnote{Graded algebras with (\ref{hnondegenerate}) are
often called non-degenerate in the literature, however, we use the
latter term in a  different sense
 (see Section~\ref{sec:multipliercross}). }. As to the first one,
$\D _g ^2 = \D _g$ is a direct consequence of the equality
\begin{equation}\label{likeprepr0}
B_g \B _{g\m } B _g = \B _g \; (\forall g\in G),\end{equation}
which is immediately seen to be true  for twisted partial crossed
products (see Section~\ref{sec:multipliercross}). The latter
equality also has several other useful consequences  given in
Lemma~\ref{easy}, and it enters  as one of the ingredients of our
criteria. In particular, it  implies that the $\D _g$'s pairwise
commute. The other ingredient deals with the $u_g$'s and $v_g$'s
we are looking for. To find them one takes the linking algebra
(also called Morita context ring)
$$\C_g =
\begin{pmatrix} \D_g   & \B_g\\ \B_{g\m} & \D_{g\m}
\end{pmatrix},$$ and it turns out that the elements $u_g$ and
$v_g$  are multipliers of $\C _g$ with $$u_g v_g = e_{11}=
\begin{pmatrix} 1 &
0 \\ 0 & 0 \end{pmatrix} \; \; \mbox{and} \; \;  v_g u_g = e_{22} =\begin{pmatrix} 0   & 0\\
0  & 1 \end{pmatrix}.$$ The matrices $e_{11}$ and $ e_{22}$  are
not elements of $\C _g $ in general, but they can be easily
interpreted as multipliers (see
Section~\ref{sec:multipliercross}).  That  crossed products by
twisted partial actions 
possess
 this property is shown in
Proposition~\ref{multipcross}.\\

The consequences of the existence of such multipliers in the
linking algebra of a given Morita context is analyzed in
Section~\ref{sec:Morita}. An isomorphism between the rings in the
 context comes out as expected, as well as a couple of maps with a
list of asso\-ci\-a\-ti\-vi\-ty properties (see
Proposition~\ref{u,v}) being used in further considerations. The
next task is to get out from the constraints of a single Morita
context in order to construct the multipliers which will form the
partial twisting. This is done in Section~\ref{sec:constructing}.
Then the criteria, Theorem~\ref{criteria}, is proved in
Section~\ref{sec:criteria} with the  restriction of homogeneous
non-degeneracy of $\B .$  The condition of the existence of the
$u_g$'s and $v_g$'s can be replaced by a more straightforward
condition, which in practice is easier to check,  if we assume
that each $\D _g$ is an $s$-unital ring (see Theorem~\ref{criteria2}).\\

The main result of \cite{E0} (Theorem 7.3) says that a second
countable $C^*$-algebraic bundle $\B $, which satisfies a certain
regularity condition, can be obtained as a twisted partial
$C^*$-crossed product. The regularity condition is based on  the
existence of a partial isometry $u$ in the multiplier algebra of a
certain linking algebra such that $u u^* = e_{11}$ and $u^* u =
e_{22}.$ If the unit fibre $\B _1$ is stable, i. e. $\B _1$ is
isometrically isomorphic to $\B _1 \otimes {\mathcal K},$ where
${\mathcal K}$ is the $C^*$-algebra of the compact operators on an
infinite dimensional separable Hilbert space, then using
\cite{BGR}, the regularity condition is guaranteed (see
\cite[Propositions 5.2, 7.1 ]{E0}). In \cite{BGR} the authors
prove  the interesting fact that two strongly Morita equivalent
$C^*$-algebras $\A $ and $\A '$ with countable approximate units
are necessarily stably isomorphic, i. e.  $\A \otimes {\mathcal K}
\cong \A ' \otimes {\mathcal K}.$ It is also shown that this fails
in the absence of countable approximate units. The condition of
$\B $ to be second countable serves to ensure the existence of
countable approximate unities (or equivalently the existence of
strictly positive elements) which is needed when using \cite{BGR}.\\

It is our next main goal  to establish an algebraic version of
Theorem 7.3 from \cite{E0}. Any $C^*$-algebra has approximate
unities, which are ``units  in limit'', and it is  an analytic
property. Our algebraic analogue of this are the local units (see
Section~\ref{sec:fdmaps} for the definitions). The property of
having countable approximate unities then corresponds to that of
having a countable set of local units. It is easily seen that such
a set of idempotents can be orthogonalized (see
Section~\ref{sec:fdmaps}), and we shall work in the more general
situation of rings with orthogonal local units. One also should
note that any closed ideal in a $C^*$-algebra is again a
$C^*$-algebra, so it possesses approximate units. For an abstract
ring with local units the analogue is not true, a two-sided ideal
in such a ring may not have local units. So we deal with a
$G$-graded $\B $ with enough local units, which means that each
$\D _g$ has orthogonal local units. Further, $\K $ can be viewed
as the $C^*$-direct (also called inductive) limit of the matrix
algebras $M_n(\Cc ),$ in which $M_n(\Cc )$ is embedded into the
left-upper corner of $M_{n+1}(\Cc ).$ The algebraic direct limit
of the $M_n(\Cc )$'s is the algebra ${\rm FMat}_{\omega }(\Cc )$
of all $\omega \times \omega$-matrices over $\Cc $ which have only
finitely many non-zero entries, where $\omega$ stands for the
first infinite ordinal. Since we do not impose restrictions on the
cardinatlities of the orthogonal units, our analogue of $\K $ will
be ${\rm FMat}_{X }(k ),$
where  $X$ is an appropriate infinite  set of indices  and $k$ is the base ring.\\

Our aim is to prove Theorem~\ref{stableteo} which says that if $\B
$ is a $G$-graded $k$-algebra with enough local units which
satisfies (\ref{likeprepr0}), then  ${\rm FMat}_{X }(\B ),$
endowed with the grading directly extended from $\B ,$ is the
crossed product by a twisted partial action of $G$ on ${\rm
FMat}_{X }(\B _1).$  Our criteria permits us to work in a single
Morita context $(\R , \R ' , \, _ {\R}M_{{\R}'},\, _ {\R '}M'_{\R
}).$ The main tool is the concept of a finitely determined map,
defined in Section~\ref{sec:constructing}. In the same section we
build up several finitely determined maps in order to use them in
the Eilenberg trick which results into a finitely determined
isomorphism $\psi : \R ^{(X)} \to \M ^{(X)}$ (see
Proposition~\ref{Eilenberg}), where $X$ is an appropriate infinite
index set.  Then in Section~\ref{sec:stab} the main step is done
in Theorem~\ref{uv-for-Morita}, in the proof of which  the map
$\psi $ is interpreted as what we call a  row and column summable
$X \times X$-matrix over ${\rm Hom} _{\R}(\R, M),$ which can be
used to define the multiplier $u$ (the analogue of the above
$u_g$) of the linking algebra  of the Morita context $({\rm
FMat}_X(\R), {\rm FMat}_X(\R '), {\rm FMat}_X(\M), {\rm
FMat}_X(\M')).$ The multiplier $v$ can be  defined by using the
inverse of $\psi $ which is also  finitely determined. In practice
we need only to give two certain maps and
Proposition~\ref{u,vexistiff} guarantees the existence of  $u$ and
 $v.$  Theorem~\ref{stableteo} is then a quick consequence of
Theorem~\ref{uv-for-Morita} thanks to the criteria. An algebraic
version of the above mentioned stable isomorphism result from
\cite{BGR} is given in Corollary~\ref{isocor}, which is an
immediate by product of Theorem~\ref{uv-for-Morita} (see also \cite{PereAra}).\\

In Section~\ref{sec:exemplos} we show by means of  examples with
uncountable local units  that Theorem~\ref{uv-for-Morita} and
Theorem~\ref{stableteo} fail if we abandon the orthogonality
condition on local units.\\



\end{section}

\begin{section}{Associativity of crossed products by twisted partial
actions}\label{sec:assoc}

In all what follows  $k$ will be a commutative associative unital
ring, which will be the base ring for our algebras. The latter
will be assumed to be associative  and non-necessarily unital. Let
 $\A$  be such an  algebra.  We
remind that  the multiplier algebra $\M(\A)$ of $\A$ is the
 set $$\M(\A)= \{(R,L) \in {\rm End}(_{\A} \A) \times {\rm End}(\A_{\A}) :
(aR)b = a(Lb) \, \mbox {for all}\, a,b
 \in \A \}$$ with component-wise addition and multiplication (see
\cite{DE}  or \cite[3.12.2]{Fillmore} for more details). Here we
use the right hand side notation   for homomorphisms of left
$\A$-modules, while for homomorphisms of right modules  the  usual
notation shall be used. In particular,  we write $a \mapsto
 a R $ and   $a \mapsto L a$  for
$R : _{\A}\A  \to _{\A}\A,$  $L : {\A}_{\A} \to {\A}_{\A}$ with $ a \in \A.$
For a multiplier $w = (R,L) \in \M(\A)$ and $a \in \A$ we  set $a w = a R$ and
$w a = L a.$ Thus
 one always has $(a w ) b = a (w b)$ $(a,b \in \A).$    The first
(resp. second) components of the elements of $\M(\A)$ are called right (resp.
 left) multipliers of $\A$. \\

\begin{defi}\label{def1} A twisted partial action of a group  $G$
 on  $\A$ is a triple
$$\Theta = (\{\D_g \}_ {g \in G}, \{{\theta}_g \}_ {g \in G},
 \{w_{g,h} \}_ {(g,h) \in G\times G}),$$ where for each $g \in G,$ ${\D}_g$ is
a two-sided ideal in $\A,$ $\theta_g$ is an isomorphism of $k$-algebras
${\D}_{g\m} \to  {\D}_g,$ and for each $(g,h) \in G\times G,$ $ w_{g,h}$ is an
invertible element from $\M(\D_g \cdot \D_{gh}),$ satisfying the following
postulates, for all $g,h$ and $t$ in $G:$\\

 (i)  $\D_g^2 = \D_g, \; \D_g \cdot \D_h = \D_h \cdot \D_g;$\\

 (ii) $\D_1 = \A$ and ${\0}_1$ is the identity map of $\A;$\\

 (iii) ${\0}_g({\D}_{g\m} \cdot {\D}_{h}) = {\D}_{g} \cdot {\D}_{gh};$\\

 (iv) ${\0}_g \circ {\0}_h (a) = w_{g,h} {\0}_{gh}(a) w\m_{g,h}, \; \forall
a \in  {\D}_{h\m} \cdot {\D}_{h\m g\m};$\\

 (v) $w_{1,g} = w_{g,1} = 1;$\\

 (vi) $\0_g(a w_{h,t}) w_{g, ht} = \0_g(a) w_{g,h} w_{gh,t}, \;
\forall a \in \D_{g\m} \cdot \D_h \cdot \D_{ht}.$

\end{defi}

 Some comments are needed concerning  the above definition.  It obviously
follows from  (i) that  a finite  product of ideals $\D_g \cdot
\D_h \cdot \ldots$ is  idempotent, and $$\0 _g (\D_{g\m} \cdot
\D_h \cdot \D_f) = \D_{g} \cdot \D_{gh} \cdot \D_{gf}$$ for all
$g,h,f \in G,$ by (iii). Thus all multipliers in (vi) are
applicable. Given an (associative) algebra $\I$ with $\I^2=\I,$ by
\cite[Prop. 2.5]{DE} one has

\begin{equation}\label{L-R} (w x) w' = w(xw')
\end{equation}

\noindent for any  $w, w' \in \M(\I), $   $x \in \I.$ This explains the absence
of brackets in the right hand side of (iv).\\

Observe also that (iii) implies

$${\0}\m _g({\D}_{g} \cdot {\D}_{h}) = {\D}_{g\m } \cdot {\D}_{g\m h},$$
for all $g,h \in G.$\\


\begin{defi}\label{crossed} Given a twisted partial action $\Theta$  of  $G$
on  $\A,$ the crossed product  $\A \stimes_{\Theta} G$ is the direct sum:
$$\bigoplus_{g \in G} \D_g \de_g,$$ in which the $\de_g$'s are symbols. The
multiplication is defined by the rule: $$(a_g {\de}_g) \cdot (b_h {\de}_h) =
{\0}_g ( {\0}\m_{g}(a_g)b_h) w_{g,h} {\de}_{gh}.$$
\end{defi}

\noindent Here $w_{g,h}$ acts  as a right multiplier on ${\0}_g (
{\0}\m_{g}(a_g)b_h) \in \0_g(\D_{g\m}\cdot \D_h) = \D_g \cdot
\D_{gh}.$\\

An element $a$ in an algebra $\A$ obviously
 determines the multiplier $(R_a, L_a) \in \M(\A)$ where $x R_a = x a$ and
 $L_a x = ax$ $(x \in \A).$  If $\I$
 is a two-sided ideal in $\A$ then this multiplier evidently restricts to
one of $\I$ which shall be denoted by the same
 pair of symbols $(R_a, L_a).$ Given an isomorphism $\af : \I \to \J$ of
algebras and a multiplier $u = (R,L)$ of $\I,$ we have that
$u^{{\af}} = ({\af}\m R \af, \af L {\af}\m)$ is a multiplier of
$\J.$ Observe that the effect of ${\af}\m R \af$ on $x \in \I$ is
$x \cdot ({\af}\m R {\af} ) = \af( {\af}\m (x) \cdot R ),$ whereas
$
(\af L {\af}\m) \cdot x = \af (L \cdot {\af}\m (x)).$\\

Before proving the associativity of  $\A \stimes_{\Theta} G$ we
separate two  technical equalities.

\begin{lema}\label{assoclemma} We have the following two properties:\\

(i) \;\;\;\;\;\;\;\; \;\;\;\;\;\;\;\;\;\;\;\;\;\;\;\;\;\;\;\;\;\;\;\;$ a
\0_h(\0\m_h(b)c) = \0_h(\0\m_h(ab)c),$\\

\noindent for any $a,c \in \A,$ $b \in \D_h$ and $h \in G;$\\

(ii)  \;\;\;\;\;\;\;\;\;\;\;\; \;\;\;\;\;\;\;\;\;\;\;
$[\0\m_{gh} (w_{g,h} \0_{gh} (x))]c =   \0\m_{gh} (w_{g,h} \0_{gh} (xc)),$\\

\noindent for any  $x \in \D_{h\m} \cdot \D_{h\m g\m},$  $g,h \in G$ and
$c \in \A.$

\end{lema}

\p (i)  Since $\0_h : \D_{h\m} \to \D_h$ is an isomorphism,
$(\0\m_h R_c{\0}_h,\0_h L_c{\0}\m_h)$ is a multiplier of $\D_h,$ and applying
(\ref{L-R})  for the multipliers  $(\0\m_h R_c{\0}_h,\0_h L_c{\0}\m_h)$ and
$(R_a, L_a)$ we see that $L_a \cdot(b \cdot (\0\m_h R_c{\0}_h)) = (L_a \cdot b)
\cdot (\0\m_h R_c{\0}_h).$  But this is precisely what (i) does say.\\

(ii) By (iii) of Definition~\ref{def1}, $\0_{gh},$ restricted to
$\D_{h\m} \cdot \D_{h\m g\m},$ gives an isomorphism $\D_{h\m}
\cdot \D_{h\m g\m} \to \D_g \cdot \D_{gh}.$ Hence $w_{g,h}^{
{\0}\m _{gh}} $ is a multiplier of $\D_{h\m} \cdot \D_{h\m g\m},$
and  combining it in (\ref{L-R}) with the multiplier $(R_c, L_c)$
we obtain  $$ [(\0\m_{gh} w_{g,h} \0_{gh}) \cdot x]\cdot R_c =
(\0\m_{gh} w_{g,h} \0_{gh}) \cdot [x \cdot R_c], $$ which is
exactly (ii). \fim \\

\begin{teo}\label{assoc} The crossed product  $\A \stimes_{\Theta} G$
is associative.
\end{teo}

\p  Obviously, $\A\stimes_{\Theta} G$ is associative if and only if

\begin{equation}\label{ght}
(a \delta_g b \delta_h ) c \delta_t = a \delta_g ( b \delta_h c \delta_t)
\end{equation}

\noindent for arbitrary $g, h, t \in G$ and $a \in \D_g, b \in \D_h,
c \in \D_t$. Computing  the left
hand side of the above equality, we have

$$ (a \delta_g b \delta_h ) c \delta_t =
\0 _g ( \0\m_{g}(a) b) w_{g,h} \delta_{gh} c  \delta_t =
 \0 _{gh} \{ \0\m_{gh}  [\0 _g ( \0\m_{g}(a) b) w_{g,h}] c \}  w_{gh,t}
\delta_{ght}.$$ On the other hand,

\begin{align*}
 & a \delta_g (b \delta_h  c \delta_t) =
a \delta_g \; \0_h(\0\m_h(b)c) w_{h,t} \delta_{ht}
=  \0_g [ \0\m_{g}(a) \0_h(\0\m_h(b)c) w_{h,t}] w_{g,ht} \de_{ght}
=\\ & \0_g [ \0\m_{g}(a) \0_h(\0\m_h(b)c) ] w_{g,h} w_{gh,t} \de_{ght},
\end{align*}

\noindent by the co-cycle equality (vi) of Definition~\ref{def1}, taking into
account that  $$ \0\m_{g}(a) \0_h(\0\m_h(b)c) \in
 \D_{g\m}\cdot \0_h(\D_{h\m}\cdot \D_t) = \D_{g\m} \cdot \D_h \cdot \D_{ht}.$$
The last equality is obtained by using (iii) of
Definition~\ref{def1}. Comparing the two sides of (\ref{ght}) we
may cancel $w_{gh,t}$ as it is invertible. Observing also that
$\0\m_g(a)$ runs over $\D_{g\m}$ when $a$ runs over $\D_g,$ we
have that (\ref{ght}) holds if and only if

\begin{equation}\label{gh}
\0 _{gh} \{ \0\m_{gh}  [\0 _g ( a b) w_{g,h}] c \} =
 \0_g [ a \0_h(\0\m_h(b)c) ] w_{g,h}
\end{equation}

\noindent is verified for any  $g, h \in G$ and $a \in \D_{g\m}, b
\in \D_h, c \in \A$ (take $t =1$ to see that $c$ may be arbitrary
in $\A$). Applying (i) of Lemma~\ref{assoclemma} to the right hand
side of (\ref{gh}), we have $\0_g [ a \0_h(\0\m_h(b)c) ] w_{g,h} =
\0_g [  \0_h(\0\m_h(ab)c) ] w_{g,h}.$ Now  $y = ab$ lies in
$\D_{g\m} \D_h,$ and therefore it is enough to show that

\begin{equation}\label{enough}
\0 _{gh} \{ \0\m_{gh}  [\0 _g (y) w_{g,h}] c \} =
 \0_g [  \0_h(\0\m_h(y)c) ] w_{g,h}
\end{equation}

\noindent is satisfied  for arbitrary  $g, h \in G,$
$y \in \D_{g\m} \cdot  \D_h, c \in \A.$ Write $x =
 \0\m_h(y) \in \0\m_h(\D_{g\m} \cdot \D_h)= \D_{h\m} \cdot \D_{h\m g\m}.$
Then by (iv) of Definition~\ref{def1} the left hand side of (\ref{enough})
becomes
$$\0 _{gh} \{ \0\m_{gh}  [\0 _g \circ \0_h(x) w_{g,h}] c \}=
 \0 _{gh} \{ \0\m_{gh}  [w_{g,h}\0 _{gh} (x) ] c \}.$$ Applying  (ii) of
Lemma~\ref{assoclemma} we see that the latter  equals to
$ \0 _{gh} \{ \0\m_{gh}  [w_{g,h}\0 _{gh} (xc) ]  \} =   w_{g,h}\0 _{gh} (xc).$ Thus
taking $z = xc,$ (\ref{enough}) becomes

$$  w_{g,h}\0 _{gh} (z) = \0_g [  \0_h(z) ] w_{g,h}$$

\noindent with arbitrary   $g, h \in G,$  $z \in \D_{h\m} \cdot
\D_{h\m g\m},$ but this is (iv) of Definition~\ref{def1} \fim \\

\end{section}

\begin{section}{Some multipliers related to
crossed products}\label{sec:multipliercross}

 We first derive   some consequences from the definition of a twisted
partial action  $\Theta.$ Taking $h=g\m, t=g$ in (vi) of Definition~\ref{def1}
we easily obtain that

\begin{equation}\label{1}
\0_g(a w_{g\m,g}) = \0_g(a)w_{g,g\m}
\end{equation} for any $g\in G, a \in \D_{g\m}.$
Applying the multiplier
$w\m_{g,g\m}$ to the both sides of the above equality and replacing
$a$ by $a w\m_{g\m,g},$
one has
\begin{equation}\label{2}
\0_g(a w\m_{g\m,g}) = \0_g(a)w\m_{g,g\m}
\end{equation} with any $g\in G, a \in \D_{g\m}.$  From (\ref{1}) and (\ref{2})  we also obtain
\begin{equation}\label{3}
\0 \m _g(a w_{g,g\m}) = \0 \m _g(a)w_{g\m,g} \; \; \forall g\in G, a
\in \D_{g};
\end{equation}  and
\begin{equation}\label{4}
\0 \m _g(a w\m_{g,g\m}) = \0 \m _g(a)w\m_{g\m,g} \; \forall g\in G, a
\in \D_{g}.
\end{equation}

Item (iv) of Definition~\ref{def1} implies $$\0_g \circ \0_{g\m} (x) =
w_{g,g\m}x w\m_{g,g\m}$$ with $x \in \D_g.$ Taking
$x = \0  _g (a), a \in \D _{g\m}$ and using (\ref{2}), one comes to
 \begin{equation}\label{1'}
\0_g(w_{g\m,g} a) = w_{g,g\m} \0_g(a) \;\; \forall g\in G, a \in \D_{g\m}.
\end{equation}

Similarly as above one has the following equalities:
\begin{equation}\label{2'}
\0_g(w\m_{g\m,g} a) = w\m_{g,g\m} \0_g(a) \;\; \forall g\in G, a \in \D_{g\m};
\end{equation}
\begin{equation}\label{3'}
\0 \m _g( w_{g,g\m} a) = w_{g\m,g}  \0 \m _g(a) \;\; \forall g\in G, a
\in \D_{g};
\end{equation}  and
\begin{equation}\label{4'}
\0 \m _g(w\m_{g,g\m } a) = w\m_{g\m,g} \0 \m _g(a) \; \;\forall g\in
G, a \in \D_{g}.
\end{equation}

Using 
the equality 
$\0_{g\m} \circ \0_{g} (x) = w_{g\m,g}x w\m_{g\m,g}$  $(x \in \D_{g\m}),$
we have

\begin{equation}\label{5}
{\0}\m _g(a) = w\m_{g\m,g} \0_{g\m}(a) w_{g\m,g} \; \; \forall g\in G,
a \in \D_{g}.
\end{equation}\\

Given a  $G$-graded algebra $\B = \bigoplus_{g\in G} \B_g,$ the
set
$$\C_g = \begin{pmatrix}
\B_g \B_{g\m}  & \B_g\\ \B_{g\m} & \B_{g\m} \B_g \end{pmatrix}$$
is evidently a $k$-algebra with respect to the usual matrix
operations. Here $\B_g \B_{g\m}$ stands for the $k$-span of the
products $a b$ with $a \in \B_g, b\in   \B_{g\m},$ and similarly
$\B_{g\m} \B_g$ is defined. Obviously $\C_g$ is contained in the
algebra of all $2\times 2$ matrices over $\B,$ and therefore, if
$\B$ is unital, the elementary matrices $e_{11}=
\begin{pmatrix} 1&0\\ 0&0
\end{pmatrix}$ and $e_{22} = \begin{pmatrix} 0&0\\ 0&1
\end{pmatrix}$ can be seen as multipliers for $\C_g.$ However even
if $\B$ is not a unital algebra, one can define the multipliers
$e_{11}$ and
$e_{22}$ of ${\C}_g$ in the natural way:\\

$$\begin{pmatrix} x & y\\ y'  & x' \end{pmatrix} \cdot
\begin{pmatrix} 1 & 0 \\ 0& 0\end{pmatrix} =
\begin{pmatrix} x & 0\\ y' &0 \end{pmatrix}, \;\;\;
\begin{pmatrix} 1 & 0 \\ 0& 0\end{pmatrix} \cdot
\begin{pmatrix} x & y\\ y' & x' \end{pmatrix}  =
\begin{pmatrix} x & y\\ 0 &0 \end{pmatrix},$$

$$\begin{pmatrix} x & y\\ y' & x' \end{pmatrix} \cdot
\begin{pmatrix} 0 &0 \\ 0& 1\end{pmatrix} =
\begin{pmatrix} 0 & y\\ 0 &x' \end{pmatrix}, \;\;\;
\begin{pmatrix} 0 & 0 \\ 0& 1 \end{pmatrix} \cdot
\begin{pmatrix} x & y\\ y' & x' \end{pmatrix}  =
\begin{pmatrix} 0 & 0\\ y' & x' \end{pmatrix},$$ where $x\in \B_g \B_{g\m},
x' \in \B_{g\m} \B_g, y \in \B_g, y' \in \B_{g\m}.$\\

In the case of a crossed product 
$\B = \A \stimes_{\Theta} G$
the algebra $\C_g$ is obviously of the form\\

$$\begin{pmatrix}
\D_g \de_1 & \D_g \de_g\\

\D_{g\m}\de_{g\m} & \D_{g\m}\de_1 \end{pmatrix}.$$\\  
If  we assume that each algebra $\D _g $ is unital  with the unity element denoted by $1_g,$ then $\C _g$ is also unital and, consequently, $\C _g$ can be identified  with its multiplier algebra. Then taking the elements 
$$ u_g=  \begin{pmatrix} 0 & 1_g \de_g\\ \;\; 0 & 0\end{pmatrix} \;\;\; \mbox{and} \;\;\;
v_g = \begin{pmatrix} 0& 0 \\ w\m_{g\m , g} \de _{g\m}& 0\end{pmatrix},$$ we see that $\C_g$ possesses multipliers $u_g $ and $v_g$ which satisfy the equalities  $u_g v_g = e_{11}$ and $v_g u_g = e_{22}.$ This however holds for general crossed products:

\begin{prop}\label{multipcross} Given a crossed product $\A \stimes_{\Theta} G$
 by a twisted partial action $\Theta,$ for each $g \in G$ there exist
multipliers
$u_g, v_g \in \M(\C_g)$ such that
$u_g v_g = e_{11}$ and $v_g u_g = e_{22}.$
\end{prop}

\p Define the multipliers $u_g$ and $v_g$ as follows.

$$ \begin{pmatrix} a\de_1 & b\de_g\\ \;\; c \de_{g\m} & d \de_{1}\end{pmatrix}
u_g = \begin{pmatrix} 0& a \de_g \\ 0&c w_{g\m,g} \de_1
\end{pmatrix},$$

$$u_g \begin{pmatrix} a\de_1 & b\de_g\\ \;\; c \de_{g\m} & d \de_{1}
\end{pmatrix} =
\begin{pmatrix} \0_g(c) w_{g,g\m} \de_1 & \0_g(d) \de_g \\ 0& 0
\end{pmatrix},$$

$$ \begin{pmatrix} a\de_1 & b\de_g\\ \;\; c \de_{g\m} & d \de_{1}\end{pmatrix}
v_g=
\begin{pmatrix} b \de_1& 0 \\ d w_{g\m,g}^{-1} \de_{g\m}&0 \end{pmatrix},$$

$$ v_g \begin{pmatrix} a\de_1 & b\de_g\\ \;\; c \de_{g\m} & d \de_{1}
\end{pmatrix} =
\begin{pmatrix} 0& 0 \\ {\0}\m _g (a w_{g,g\m}^{-1}) \de_{g\m}& {\0}\m
_g(b)\de_1
\end{pmatrix}.$$

Check first that $u_g$ is a multiplier. Write
$$x = \begin{pmatrix} a\de_1 & b\de_g\\ \;\; c \de_{g\m} & d \de_{1}
\end{pmatrix}
\;\;\; \mbox{and} \;\;\; y = \begin{pmatrix} a'\de_1 & b'\de_g\\
\;\; c' \de_{g\m} & d' \de_{1}
\end{pmatrix}$$  By (\ref{1}) we easily have

\begin{align*} &(x y ) u_g = \begin{pmatrix} 0& (a a'+
b \0_{g}(c')w_{g,g\m})\de_g\\
0 & (c \0_{g\m}(a') + d c')w_{g\m,g} \de_1 \end{pmatrix} = x(y
u_g).
\end{align*} Furthermore, for the $(1,1)$-entry of $u_g (xy),$
by (iv) of Definition~\ref{def1}, one has
\begin{align*}  & \0_g( c ) \; \0_g \circ \0_{g\m} (a') w_{g,g\m} +
\0_g(d c')w_{g,g\m} =
              \\& \0_g( c ) w_{g, g\m} a' w \m _{g,g\m} w_{g,g\m} +
\0_g(d )\0_g(c' )w_{g,g\m},
\end{align*} which is the $(1,1)$-entry of  $(u_g x)y.$ On the other hand,
 by (\ref{1}) the $(1,2)$-entry of $u_g (xy)$ equals
\begin{align*}  &(\0_g( c ) \; \0_g\circ \0_{g\m} (b') w_{g,g\m} +
\0_g(d d' ))\de_g =\\
                 & (\0_g( c ) w_{g,g\m} b'
+ \0_g(d)\0_g(d'))\de_g, \end{align*} using again (iv) of
Definition~\ref{def1}. It is easily seen that this is the
$(1,2)$-entry of $(u_g x)y,$ implying $u_g (xy) =(u_g x)y.$ One
also has

\begin{align*}  x(u_g y) =
\begin{pmatrix} a \0 _g (c') w_{g,g\m} \de _1 & a \0_g( d' ) \de_g\\
c\0_{g\m}(\0 _g (c') w_{g,g\m})\de _{g\m} & c \0 _{g\m} ( \0 _g
(d')) w_{g\m,g}\de _1
\end{pmatrix} =
(x u_g) y
\end{align*} applying  once more (\ref{1}) and (iv) of Definition~\ref{def1}. This shows that $u_g$
is a multiplier of ${\C}_g.$\\

Similarly, using (\ref{2}) one has $(xy)v_g = x(yv_g),$ and by
(\ref{4}) and (\ref{5}) it is easily seen that $(xv_g )y = x (v_g
y).$ As to the equality $v_g (xy) = (v_g x)y,$  the $(2,1)$-entry
of $v_g (xy)$ is

\begin{align*}
& {\0}\m_g(a a' + b \0_g (c' ) w_{g,g\m})w\m_{g\m,g} \de_{g\m}=\\
& ({\0}\m _g(a){\0}\m _g(a') w\m_{g\m,g} + {\0}\m_g(b)c')\de_{g\m},
\end{align*} by (\ref{4}). Applying (\ref{5}) this equals to the $(2,1)$-entry of $(v_g x)y.$
Using (\ref{4}) and (\ref{5}) we also see that the the
$(2,2)$-entries coincide, which shows that $v_g (xy) = (v_g x)y.$
Thus $v_g \in \M({\C}_g).$
By easy calculations we have  $u_g v_g = e_{11}$ and $v_g u_g = e_{22}.$ \fim \\

 We need one more property of a crossed product  $\A \stimes_{\Theta} G$
by a twisted partial action, which is in fact an  immediate consequence of its
definition. More precisely,  we easily see that
$\D_g \de_g \cdot \D_{g\m} \de_{g\m} \cdot \D_g \de_g =
(\D_g^2) w_{g,g\m} \de_e \cdot \D_g \de_g = \D_g \de_e \cdot \D_g \de_g =
 \D_g \de_g,$ as $\D_g$ is idempotent and
$w_{g,g\m}$ is an invertible multiplier of $\D_g.$ Thus the
equality

\begin{equation}\label{likeprepr1}  \D_g \de_g \cdot \D_{g\m} \de_{g\m}
\cdot \D_g \de_g
= \D_g \de_g
\end{equation} holds for any $g \in G.$\\

We shall work with a mild restriction on $\Theta,$ supposing that
each domain $\D_g$ is left and right non-degenerate. Following
\cite{DE} and \cite{GS} we say that an algebra $\I$
 is   {\em right non-degenerate} if for any $0 \neq x \in \I$ one has
$x \I \neq \{0\},$ and similarly is defined the concept of a {\em left
non-degenerate} algebra. If each $\D_g$ is left and right non-degenerate then
we evidently have

\begin{equation*} x \neq 0 \Longrightarrow x\cdot \D_{g\m} \de_{g\m} \neq 0
\;\; \mbox{and}\;\;  \D_{g\m} \de_{g\m} \cdot x \neq 0 \;\;
\forall x \in \D_g \de_{g}.
\end{equation*} For a general $G$-graded algebra $\B =
\bigoplus_{g \in G}\B_g$
this corresponds  to the property

\begin{equation}\label{hnondegenerate} x\neq 0 \Longrightarrow
x \cdot \B_{g\m}  \neq 0 \;\; \mbox{and}\;\; \B_{g\m} \cdot x \neq 0 \;\;
\forall  x \in \B_g.
\end{equation}

We shall say that a $G$-graded algebra $\B$ is {\em homogeneously non-degenerate} if $\B$ satisfies
(\ref{hnondegenerate}). It turns out that among the  homogeneously non-degenerate graded algebras the crossed
products by twisted partial actions are distinguished by the property (\ref{likeprepr1}) and the existence of
the multipliers like in Proposition~\ref{multipcross}.  The precise statement will
be given in Theorem~\ref{criteria}.\\

 Given a graded algebra $\B,$ write $\D_g = \B_g \B_{g\m}.$ Observe
that $\B_g$ is a $(\D_g,\D_{g\m})$-bimodule, whereas $\B_{g\m}$ is
a $(\D_{g\m},\D_{g})$-bimodule, so that $(\D_g,\D_{g\m},  \B_g,
\B_{g\m})$ together with the surjective maps $\B_{g}\otimes
\B_{g\m} \to \D_{g}$ and $\B_{g\m} \otimes \B_{g} \to \D_{g\m},$
determined by the   multiplication in $\B,$ form a Morita context. In general
 we know nothing about the
existence of an identity element in $\D_g,$ however in the
situation in which we are interested in, the $\D_g$'s are
idempotent rings and the bimodules ${\B}_g, {\B}_{g\m}$ are
unital, and so by a result of \cite{GS}, $\D_{g\m}$ and $\D_g$ are
Morita equivalent (see Section~\ref{sec:fdmaps} for the precise
definitions). The idempotence of the $\D_g$'s,   as well as the
property of the ${\B}_g$'s being unital modules, is an immediate
consequence of the analogue of (\ref{likeprepr1}) for graded
algebras:

\begin{equation}\label{likeprepr2}  \B_g  \cdot \B_{g\m}  \cdot \B_g
= \B_g  \;\; \forall g \in G.
\end{equation}\\

  The existence of multipliers of $\B$, like in Proposition~\ref{multipcross},
turns out to be a powerful tool: we shall have  that $\D_{g\m}$
and $\D_g$  are isomorphic and, the obtained isomorphisms will fit
into a twisted partial action $\Theta$ of $G$ on $\B_1$ so that
$\B \cong  \B_1 \stimes_{\Theta} G.$ In order to carry out this
idea we need some technical tools, which we develop in the next
sections.\\

 \end{section}

\begin{section}{Morita context and multipliers}\label{sec:Morita}

Let $(\R, {\R}', M, M', \tau, {\tau}')$
be a Morita context in which $\R$ and ${\R}'$ are some (non-necessarily unital)
algebras. We recall that this means that  $M$ is an $\R$-$\R'$-bimodule, $M'$
is an $\R'$-$\R$-bimodule, $\tau : M\otimes_{\R'} M' \to \R$ is an
$\R$-bimodule
map, $\tau' : M'\otimes_{\R} M \to \R'$ is an $\R'$-bimodule
map,  such that $$
  \tau(m_1\otimes m')\;m_2 = m_1\; \tau'(m'\otimes m_2),\qquad
\forall m_1,m_2\in
M,\; m'\in M',$$
  and
  $$\tau'(m_1'\otimes m)\;m'_2 = m'_1\; \tau(m\otimes m'_2),\qquad \forall
m'_1,m'_2\in M',\; m\in M.$$ One can construct the {\em linking algebra}
of the Morita context, which is the set $$\C = \begin{pmatrix}
\R  & M \\ M'& \R' \end{pmatrix},$$ with the obvious addition of matrices
and multiplication by scalars, and
the matrix  multiplication
determined by the bimodule structures on $M$ and $M'$ and the maps $\tau$ and
$\tau',$ so that $m \cdot m' = \tau (m\otimes m')$ and
$m' \cdot m = \tau' (m' \otimes m).$

If we have \begin{align}\label{idempotentC} & {\R}^2 = \R, ({\R'})^2 = \R',
\R M + M\R' = M \;\;\mbox{and}\;\; \R' M' + M' \R = M',
\end{align} it is easily verified that ${\C}^2 = \C.$ Thus by
\cite[Prop. 2.5]{DE}

\begin{equation}\label{commuting-on-C}
(wx)w' = w(xw') \;\;\;\forall  w,w' \in \M(\C), x\in \C.
\end{equation}

In the case of a graded algebra $\B$ with
$\R = \D_g, \R' = D_{g\m}, M= \B_g$ and $M' = \B_{g\m},$ the linking algebra
$\C$ evidently coincides with the algebra $\C_g$ considered in the previous
section. If the graded algebra $\B$ satisfies (\ref{likeprepr2}) then
obviously (\ref{idempotentC}) is also satisfied.

Similarly, as it was done in the previous section for $\C_g,$ one
defines the multipliers $e_{11}$ and $e_{22}$ of $\C$ in the
natural way.\\
%

We are going to examine the situation in which the linking algebra $\C$
of a Morita context contains multipliers with properties announced in
Proposition~\ref{multipcross}. More precisely, suppose that we have a Morita
context  $(\R, {\R}', M, M', \tau, {\tau}')$ such that the linking algebra
$\C$ is idempotent. Suppose furthermore that there exists $u, v \in \M(\C)$
such that $uv = e_{11}$ and $vu = e_{22}.$ One may assume that
$e_{11} u= u e_{22} =u$ and $e_{22} v=v e_{11} = v.$ For it is enough to
replace $u$ by $e_{11}ue_{22}$ and $v$ by $e_{22}ve_{11},$ as
$(e_{11}ue_{22})(e_{22}ve_{11})= (uv)^4 =e_{11}$ and similarly
$ (e_{22}ve_{11})(e_{11}ue_{22})= e_{22}.$ Let $x=
\begin{pmatrix}r&m\\m'&r'\end{pmatrix}$ be arbitrary in $\C.$ Then
$$ux = e_{11}u(e_{22}x) = e_{11} u \begin{pmatrix}0&0\\m'&r'\end{pmatrix} =
\begin{pmatrix}\ast&\ast\\0&0\end{pmatrix}.$$ Moreover, by
(\ref{commuting-on-C}),
\begin{align*} &u \begin{pmatrix} 0&0 \\m'&0\end{pmatrix} =
u (xe_{11})=  (ux)e_{11} =
\begin{pmatrix} \ast&0\\0&0 \end{pmatrix},
\end{align*} and the $(1,1)$-entry of the latter matrix depends only on $m'.$ Denote it by
$u\cdot m'.$ Thus we have a map of $k$-modules
$M'\ni m' \mapsto u\cdot m' \in \R$ which we denote by $L_u.$ Similarly,
\begin{align*} &u \begin{pmatrix}0&0\\0&r' \end{pmatrix} = u (xe_{22})=
(ux)e_{22} = \begin{pmatrix}0& u\cdot r' \\0&0\end{pmatrix},
\end{align*} and we have a $k$-linear map $\R' \ni r' \mapsto u\cdot r' \in M,$
which we denote by the same symbol $L_u.$ Thus we have
\begin{equation*} u \begin{pmatrix}0&0\\m'&r'\end{pmatrix} =
\begin{pmatrix} um'& ur'\\0&0\end{pmatrix}.
\end{equation*}
Analogously,  $u,$ as a right multiplier, determines  a pair of maps
$$R_u: \R \ni r \mapsto r\cdot u \in M\;\;\mbox{and}\;\;
R_u: \M' \ni m' \mapsto m'\cdot u \in \R'$$ so that
\begin{equation*}  \begin{pmatrix}r&0\\m'&0\end{pmatrix} u =
\begin{pmatrix} 0& r u\\0&m' u\end{pmatrix}.
\end{equation*}

In a similar way the multiplier $v$ gives rise to the pairs of $k$-linear maps:
$$L_v: \R \to M', M \to \R' \;(r \mt v\cdot r, m \mt v\cdot m)$$ and
$$R_v:M\to \R, \R' \to M' \; (m\mt m\cdot v, r' \mt r' \cdot v),$$ so that

\begin{equation*} v \begin{pmatrix}r&m\\0&0\end{pmatrix}  =
\begin{pmatrix} 0& 0\\v r&v m\end{pmatrix},
\end{equation*} and
\begin{equation*}  \begin{pmatrix}0&m\\0&r'\end{pmatrix} v =
\begin{pmatrix} mv& 0\\r' v&0\end{pmatrix}.
\end{equation*}\\

\begin{lema}\label{u,v}  Let $(\R, {\R}', M, M', \tau, {\tau}')$ be a Morita
context  such that the linking algebra $\C$ is idempotent. Suppose
that the multiplier algebra $\M (\C)$ contains elements $u$ and
$v$ with $u v = e_{11}$ and $v u = e_{22}.$ Then the multipliers
$u$ and $v$ can be chosen such that $u$ determines the maps
$$ L_u: M' \ni m'  \mapsto u \cdot m' \in \R, \;\;\; \R' \ni r' \mapsto
u \cdot r'  \in M,$$
$$R_u: \R \ni r \mapsto r \cdot u \in M, \;\;\; M' \ni m' \mapsto
m'\cdot u \in \R',$$
which are isomorphisms of $k$-modules whose inverses are respectively the maps

$$L_v: \R \ni r  \mapsto v \cdot r \in M', \;\;\;
M \ni m \mapsto v \cdot m  \in \R',$$
$$R_v: M \ni m \mapsto m \cdot v \in \R, \;\;\;
\R' \ni r' \mapsto r'\cdot v \in M',$$
determined by  $v.$ Moreover, the following properties hold  with arbitrary
$r, r_1,r_2\in \R, r',r'_1,r'_2\in \R',m, m_1, m_2\in M,m',m'_1,m'_2 \in M':$
\begin{eqnarray*}
u(m' r) =(u m')r,          &(r_1 r_2)u = r_1 (r_2 u),  &(r u) m' = r(u m'),\\
u(m' m) = (u m')m,         &(m m')u = m(m' u),       &(r u)r'=r(u r'),\\
u(r'm') = (u r') m',      &(r'm')u = r'(m'u),     &(m_1'u)m_2'=m'_1(um'_2),\\
u(r'_1r'_2)= (ur'_1)r'_2, &(m'r)u  = m'(r u), &(m'u)r'=m'(u r'),
\end{eqnarray*} and the corresponding properties   involving $v:$
\begin{eqnarray*}
v(r_1r_2) =(vr_1)r_2,&(r m)v = r (m v),                 &(m v) r = m(v r),\\
v(m m') = (v m)m',    &(m r')v = m(r' v),               &(m_1v)m_2=m_1(vm_2),\\
v(r m) = (v r) m,     &(m'm)v = m'(m v),               &(r'v)m=r'(v m),\\
v(m r')= (v m)r',   &(r'_1r'_2)v =r'_1(r'_2v), &(r'v)r=r'(v r).
\end{eqnarray*}
\end{lema}
\p The pairs of maps $L_u, R_u, L_v, R_v$ are already defined.
It follows immediately from the equalities $uv=e_{11}, vu=e_{22}$ that
all these maps are isomorphisms of $k$-modules with inverses as announced
in the lemma.

Using the fact that $u$ is a left multiplier, we have
\begin{align*} & \left(u \begin{pmatrix} 0&0\\m'&0\end{pmatrix}\right)
\begin{pmatrix} r&m\\0&0\end{pmatrix} =
\begin{pmatrix} (um') r&(um')m\\0&0\end{pmatrix} = \\
& u \left(\begin{pmatrix} 0&0\\m'&0\end{pmatrix}
\begin{pmatrix} r&m\\0&0\end{pmatrix}\right) =
 \begin{pmatrix} u(m'r)&u(m m')\\0&0\end{pmatrix},
\end{align*} which gives the first two properties enounced in the lemma.
The next two properties of the first column  come from
\begin{align*} &\left(u \begin{pmatrix} 0&0\\0&r'_1\end{pmatrix}\right)
\begin{pmatrix} 0&0\\m'&r'_2\end{pmatrix} =
\begin{pmatrix} (ur'_1)m' &(ur'_1)r'_2\\0&0\end{pmatrix} =\\
& u \left(\begin{pmatrix} 0&0\\0&r'_1\end{pmatrix}
\begin{pmatrix} 0&0\\m'&r'_2\end{pmatrix}\right) =
\begin{pmatrix} u(r'm')&u(r'_1 r'_2)\\0&0\end{pmatrix}.
\end{align*}
The second column of properties involving $u$ is obtained from the
fact that $u$ is a right multiplier, whereas the last column comes
from the equality $(xu)y=x(uy)$ with appropriate $x,y \in \C.$ The
analogous properties involving $v$ are obtained similarly.\fim

\begin{remark}\label{allassoc} The list of equalities  in the  lemma above
say that  whenever we form a product of three elements one of which is $u$ or
$v$ and the other two belong to $M \cup M' \cup \R \cup \R',$ it is
associative, i.e. we may put the  brackets  arbitrarily.
\end{remark}

\begin{remark}\label{uv(iff)} It is directly seen that two pairs of some  $k$-linear maps
$$ L_u: M' \ni m'  \mapsto u \cdot m' \in \R, \;\;\; \R' \ni r' \mapsto
u \cdot r'  \in M,$$ $$R_u: \R \ni r \mapsto r \cdot u \in M,
\;\;\; M' \ni m' \mapsto m'\cdot u \in \R'$$ determine a
multiplier $u$ of $\M (\C )$ exactly when the equalities of the
three columns involving $u$ in Lemma~\ref{u,v} hold (and similarly
for $v$).
\end{remark}


\begin{remark}\label{assumptions} One can directly check that in 
Lemma~\ref{u,v} the assumption that $\C$ is idempotent can be replaced by that of $({\R}')^2 ={\R}', {\R}' M' = M' \R = M'$ or by the assumption ${\R}^2 = \R , \R M =  M{\R}' =M.$
\end{remark}

The pairs of maps $L_u$ and $R_u$ given in Lemma~\ref{u,v} permit
to establish an isomorphism between $\R'$ and $\R:$

\begin{prop}\label{iso} Let $(\R, {\R}', M, M', \tau, {\tau}')$ be a Morita
context with idempotent linking algebra $C$  such that $\M (\C)$ contains
elements $u$ and $v$ with  $u v = e_{11}$ and
$v u = e_{22}.$ Then  $u$ and $v$ can be chosen such  that the map
$$\0 : \R' \ni r' \mt ur'v \in \R$$ is an
isomorphism $k$-algebras, whose inverse is  $$\R \ni r \mt vru \in \R'.$$
\end{prop}
\p Observe first that by (\ref{commuting-on-C}) we have

\begin{align*}  \left(u \begin{pmatrix} 0&0\\0 & r'\end{pmatrix}\right)v=
u\left (\begin{pmatrix} 0&0\\0&r'\end{pmatrix} v \right) \;\; \mbox{and}\;\;
 \left(v \begin{pmatrix} r&0\\0 & 0\end{pmatrix}\right)u=
v\left( \begin{pmatrix} r&0\\0&0\end{pmatrix} u\right) ,
\end{align*} which  implies $(ur')v = u(r'v)$ and $(vr)u = v(ru)$ so that it makes sense to omit
brackets
in $ur'v$ and $vru.$

Clearly, $(u(\R'))v = (M)v= \R$ so $\0$ maps $\R'$ onto $\R.$ Moreover,
because $L_u$ and $R_v$ are $k$-isomorphisms, it follows that $\0$ is also a
$k$-isomorphism. Clearly $\R \ni x \mt vxu \in \R' $ is the inverse of $\0 .$  Thus it
remains to show that $\0$ preserves multiplication. This follows from the next calculations:
\begin{align*} & \begin{pmatrix} \0 (r'_1r'_2)&0\\0 & 0\end{pmatrix}=
u\begin{pmatrix} 0&0\\0&r'_1  \end{pmatrix}  \begin{pmatrix} 0&0\\0& r'_2 \end{pmatrix} v =
\left( u\begin{pmatrix} 0&0\\0&r'_1  \end{pmatrix}  e_{22}\right)  \begin{pmatrix} 0&0\\0& r'_2
\end{pmatrix} v =\\
&\left( u\begin{pmatrix} 0&0\\0&r'_1  \end{pmatrix} vu\right)  \begin{pmatrix} 0&0\\0& r'_2
\end{pmatrix} v =
\begin{pmatrix} \0 (r'_1)&0\\0 & 0\end{pmatrix} \begin{pmatrix} \0 (r'_2)&0\\0 & 0\end{pmatrix}.
\end{align*}
 \fim \\


Given a $G$-graded algebra $\B = \bigoplus_{g \in G} \B_g$
which satisfies (\ref{likeprepr2})  and such that for any $g \in G$
there exist $u_g, v_g \in \M (\C_g)$  with  $u_g v_g = e_{11}$ and
$v_g u_g = e_{22},$ applying Lemma~\ref{u,v}
to the Morita context $(\D_g,\D_{g\m},  \B_g, \B_{g\m}, \cdot, \cdot),$ we
have the isomorphisms  of $k$-modules
$$L_{u_g}: \B_{g\m} \to \D_g, \D_{g\m} \to \B_g, \;\;
R_{u_g}: \D_{g} \to \B_g, \B_{g\m} \to \D_{g\m}$$ with inverses
$$L_{v_g}: \D_{g} \to \B_{g\m}, \B_{g} \to \D_{g\m},
R_{v_g}: \B_{g} \to \D_g, \D_{g\m} \to \B_{g\m},$$ respectively. Moreover,
by Proposition~\ref{iso}
\begin{equation}\label{isos}
{\0}_{g} : \D_{g\m} \ni x \mt u_gxv_g \in \D_g
\end{equation}  is an isomorphism of $k$-algebras.\\

\end{section}

\begin{section}{Constructing $w_{g,h}$}\label{sec:constructing}

\begin{lema}\label{extending} Let $\B$ be a homogeneously
non-degenerate $G$-graded algebra  satisfying (\ref{likeprepr2})  and such
that for any $g \in G$ the multiplier algebra $\M (\C_g)$ contains
elements $u_g$ and $v_g$ with  $u_g v_g = e_{11}$ and
$v_g u_g = e_{22}.$
Then there are isomorphisms  of right $\B$-modules
$$\tilde L_{u_g}: \B_{g\m} \B \ni xy \mt (u_g x)y \in \B_g \B \;\;\;
(x\in \B_{g\m}, y \in \B),$$
$$\tilde L_{v_g}: \B_g \B \ni xy \mt (v_g x)y \in \B_{g\m} \B \;\;\;
(x \in \B_g, y \in \B)$$ and isomorphisms of left $\B$-modules
 $$\tilde R_{u_g}: \B \B_{g\m} \ni xy \mt x(y u_g) \in \B \B_g \;\;\;
(x \in \B, y \in \B_{g\m}),$$
$$\tilde R_{v_g}: \B \B_{g} \ni xy \mt x(y v_g) \in   \B  \B_{g\m}
\;\;\; (x \in \B, y \in \B_{g}),$$  such that
$\tilde L_{u_g}$ extends the pair of maps
$$L_{u_g}:\B_{g\m} \to \D_g, \D_{g\m} \to \B_g,$$
$\tilde R_{u_g}$ extends
$$R_{u_g}: \D_{g} \to \B_g, \B_{g\m} \to \D_{g\m},$$ and similarly for
$\tilde L_{v_g}$ and $\tilde R_{v_g}.$  Moreover, $\tilde L_{v_g}
=\tilde L\m_{u_g}$ and $\tilde R_{v_g} = \tilde R\m_{u_g}.$
\end{lema}

\p  Observe first that
\begin{equation}\label{Lg1} u_g(x y) = u_g(x) y \;\;\;  \forall x\in \B_{g\m},
y\in \B_e.
\end{equation} Indeed, by (\ref{likeprepr2}), $\B_{g\m} \D_g =
\B_{g\m} \B_{g}\B_{g\m} = \B_{g\m},$ and thus we can write
$x = \sum x'_i x''_i$ with $x'_i \in \B_{g\m}, x''_i \in \D_g.$ Then using the
first equality in Lemma~\ref{u,v}, and keeping in mind that
$x''_i y \in \D_g \B_e \subseteq \D_g,$ we have
$$ u_g(xy) = \sum u_g(x'_i x''_iy) =
\sum u_g(x'_i)x''_iy =\sum u_g(x'_ix''_i)y = u_g(x)y,$$ as desired.

We define  $\tilde L_{u_g}$ first on $\B_{g\m} \B_h$ with arbitrary
$h \in G$ by
$$\tilde L_{u_g} (xy) = (u_g x)y \;\;\;
\forall x \in \B_{g\m}, y\in \B_h,$$ and extending it additively to
arbitrary elements of $\B_{g\m} \B_h.$   Then $$\tilde L_{u_g}:\B_{g\m} \B_h
\to \D_g \B_h$$  is a well-defined map. For suppose $\sum x_i y_i =0$ for some
$x_i \in \B_{g\m}, y_i \in \B_h.$ Then for any $z \in \B_{h\m}$ applying
(\ref{Lg1}) we have
$$  (\tilde L_{u_g} (\sum x_i y_i))z =  \sum u_g(x_i)y_iz =
\sum u_g(x_iy_iz) = u_g(\sum x_iy_iz)=0,$$ as $y_iz \in \B_e.$
Since $z \in  \B_{h\m}$ is arbitrary, it follows that $(\tilde
L_{u_g}(\sum x_i y_i)) \B_{h\m} =0,$ which implies $\tilde L_{u_g}
(\sum x_i y_i) =0,$ as $\B$ is homogeneously non-degenerate.

Since $\B = \bigoplus_{h \in G} \B_h$ we can extend  $\tilde
L_{u_g}$  additively to a well-defined map $$\tilde L_{u_g} :
\B_{g\m} \B \to \D_g \B,$$ which sends $xy \mt (u_g x)y$ with $x
\in \B_{g\m}, y \in \B.$  Observe now that
$$\B_g \B_h = \B_g \B_{g\m} \B_g \B_h \subseteq \D_g \B,$$ in view of
(\ref{likeprepr2}). Consequently, $\B_g \B \subseteq \D_g \B,$ and since
obviously $\D_g \B \subseteq \B_g \B$ , one has
$\D_g \B = \B_g \B.$

Thus we obtain a map $$\tilde L_{u_g}: \B_{g\m}\B \to \B_g \B ,$$ which is a
homomorphism of right $\B$-modules by its construction. Its restriction to
$\B_{g\m}$ clearly coincides with $L_{u_g}.$ Now let
$z \in \D_{g\m}=\B_{g\m} \B_g$ be arbitrary and write
$z = \sum x_iy_i,\; x_i \in \B_{g\m}, y_i \in \B_g.$
Then using the second equality of the first column in Lemma~\ref{u,v}, we
have
\begin{align*} \tilde L_{u_g}(z) = \sum ({u_g} x_i) y_i=
\sum u_g (x_iy_i) = L_{u_g} (z),
\end{align*} so that $\tilde L_{u_g},$ when restricted to $\D_{g\m},$
coincides with $L_{u_g}.$

We proceed similarly with $R_{u_g}.$ Using the equality
$\B_{g\m} = \B_{g\m} \B_g \B_{g\m},$ we readily obtain
\begin{equation}\label{Rg1}
(xy)u_g = x(yu_g)\;\;\; \forall x \in \B_e, y \in \B_{g\m}.
\end{equation} Then for an arbitrary fixed  $h \in G$ set
$$\tilde R_{u_g} (x y) = x (y u_g) \;\;\;
\forall  x\in \B_h, y \in \B_{g\m},$$ which gives a well-defined map
$$\tilde R_{u_g} : \B_h \B_{g\m} \to \B_h \D_{g\m}.$$ Indeed, if
$\sum x_i y_i =0$ with $x_i \in \B_h, y_i \in \B_{g\m}$ then $z
x_i \in \B_1$ for any $z \in \B_{h\m}$ and using (\ref{Rg1}) we
have
$$z \sum x_i (y_i u_g) = \sum (z x_i y_i) u_g =0.$$ As above this yields
$\tilde R_{u_g} (\sum x_i y_i) =0,$ as $\B$ is homogeneously
non-degenerate.

Extending additively we come to a well-defined homomorphism of left
$\B$-modules
$$\tilde R_{u_g} : \B \B_{g\m} \to \B \D_{g\m}= \B \B_{g}.$$

The maps $\tilde L_{v_g}$ and
$\tilde R_{v_g}$ are obtained by analogous arguments. Finally,
for $x \in \B_{g\m}, y \in \B$  we obviously  have
$$\tilde L_{v_g} \circ \tilde L_{u_g} (xy) =
\tilde L_{v_g}((u_g x)y) =  (v_g u_g x)y = xy,$$ as $L_{v_g} =
L\m_{u_g}.$ It is also evident that $\tilde L_{u_g} \circ \tilde
L_{v_g}$ is identity on $\B_g \B.$ Hence $\tilde L_{v_g} = \tilde
L\m_{u_g}$ and similarly $\tilde R_{v_g} = \tilde R\m_{u_g}.$ \fim \\

We shall write $u_g x = \tilde L_{u_g} (x),  x u_g = \tilde
R_{u_g} (x), v_g x = \tilde L_{v_g} (x)$ and $x v_g = \tilde
R_{v_g} (x),$ provided that $x$ belongs to the domain of the
considered map.\\

\begin{remark}\label{convenient} Suppose that  $\B$ is  as in the above lemma. It will be
convenient to observe the following equalities which easily follow
from  Lemma~\ref{extending}:\\
\begin{align*}
& u_g(xy) = (u_g x) y \;\;\;\;\; \forall x\in \D_{g\m}, y \in \B,\\
& (xy)u_g = x(y u_g) \;\;\;\;\; \forall x\in \B, y \in \D_g,\\
& v_g (xy) = (v_g x)y \;\;\;\;\; \forall x \in \D_g, y \in \B,\\
& (xy)v_g = x(y v_g) \;\;\;\;\; \forall x \in \B, y \in
\D_{g\m}.\\
\end{align*}
\end{remark}

We shall use the  maps obtained in Lemma~\ref{extending} to
construct the multipliers $w_{g,h}.$ We list first some easy
consequences of the property (\ref{likeprepr2}).\\

\begin{lema}\label{easy} Suppose that  $\B = \bigoplus_{g \in G} \B_g$
is a graded algebra satisfying (\ref{likeprepr2}). Then for all $g,h \in G$
we have:\\

(i) $\B_{g\m} \B_g \B_h = \B_{g\m} \B_{gh},$  $\B_{g} \B_h \B_{h\m} = \B_{gh} \B_{h\m},$\\

(ii) $\D_g \D_h = \D_h \D_g,$\\

(iii) $\D_g \B_g = \B_g, \B_g \D_{g\m} = \B_g.$\\
\end{lema}
\p Clearly,
\begin{align*}&\B_{g\m} \B_g \B_h \subseteq \B_{g\m} \B_{gh} =
\B_{g\m} \B_g \B_{g\m} \B_{gh} \subseteq \B_{g\m} \B_g \B_h,
\end{align*} which gives the first equality of  (i). The second one follows similarly. Using (i)
we have
\begin{align*}&\D_g \D_h = \B_g \B_{g\m} \B_h \B_{h\m} =
\B_g \B_{g\m h} \B_{h\m} =
 \B_g \B_{g\m h} \B_{h\m g} \B_{g\m h} \B_{h\m} = \\
& \B_{h} \B_{h\m g} \B_{g\m}=
\B_h \B_{h\m} \B_g \B_{g\m}
\end{align*} and (ii) follows. Finally, item (iii) is immediate.
\fim \\

\begin{prop}\label{w} Suppose that  $\B$ is a $G$-graded algebra as in Lemma~\ref{extending}.
Then for
any $g, h \in G$ the equality $$w_{g,h} = u_g u_h v_{gh}$$ defines an
invertible multiplier of $\D_g \D_{gh}$ with $$w\m_{g,h} = u_{gh} v_h v_g.$$
\end{prop}
\p  Take arbitrary $x,y \in \D_g \D_{gh} = \D_{gh} \D_g$ and write
$x = \sum x' x''$ with $x' \in \D_{gh}$ and $x'' \in \D_g$
(we omit the indices for simplicity). Then by
Remark~\ref{convenient}
\begin{align*}  &v_{gh} (xy) =  \sum v_{gh} (x' x'' y)= \sum v_{gh} (x') x'' y=
\sum v_{gh} (x' x'')y = (v_{gh} x)y.
\end{align*} We also have
$$v_{gh}(x) \in v_{gh} (\D_{gh} \D_g) = v_{gh}(\D_{gh})\D_g  =
\B_{h\m g\m}\D_g = \B_{h\m} \B_{g\m} \D_g = \B_{h\m} \B_{g\m},$$
by Remark~\ref{convenient} and (\ref{likeprepr2}). So we can write
$v_{gh}(x) = \sum  x'  x''$ with $ x' \in \B_{h\m},  x'' \in
\B_{g\m}.$ Then  by  Lemma~\ref{extending}
\begin{align*} &u_h v_{gh} (x y) = \sum u_h ( x'  x'' y) =
\sum u_h ( x')  x'' y = \sum u_h ( x'  x'')y = (u_hv_{gh}(x))y.
\end{align*} We see that
\begin{align*}&u_h v_{gh}(x) \in u_h(\B_{h\m} \B_{g\m})= u_h(\B_{h\m})\B_{g\m} =
 \D_{h}\B_{g\m} = \D_h \D_{g\m} \B_{g\m} = \D_{g\m} \D_h \B_{g\m},
\end{align*} and we can write $u_hv_{gh}(x) = \sum x' x''$ with $x' \in
\D_{g\m}$ and $x'' \in \B.$ Then
\begin{align*} & u_g ((u_h v_{gh}x)y) = \sum u_g (x' x'' y) =
\sum u_g(x')x''y=
\sum u_g(x' x'')y = (u_gu_hv_{gh}(x))y.
\end{align*} This shows that $$u_g u_h v_{gh} (xy) =
u_g u_h v_{gh}(x)y \;\;\; \forall x,y \in \D_g \D_{gh}.$$  We also have
\begin{align*}& u_gu_hv_{gh}(\D_{gh}\D_g)=u_g(\D_{g\m} \D_h \B_{g\m})=
u_g(\D_{g\m})\D_h \B_{g\m} = \B_g \D_h \B_{g\m}\\& =
\B_g \D_h \B_{g\m}\D_{g}= \B_g \B_h \B_{h\m} \B_{g\m} \D_g=
\B_{gh}\B_{h\m g\m}\D_g =\D_{gh} \D_g.
\end{align*}

Using similar calculations one also shows that
$$(xy) u_g u_h v_{gh}  =
x(y u_g u_h v_{gh}) \;\;\; \forall x,y \in \D_g \D_{gh}$$ and
$$ (\D_{gh}\D_g) u_gu_hv_{gh}=\D_{gh} \D_g.$$

On order to prove  that $w_{g,h}= u_g u_h v_{gh}$ is a multiplier of
$\D_g \D_{gh}$ it remains to show that
$$ (x w_{g,h})y = x (w_{g,h} y) \;\;\;\;\; \forall x,y \in \D_g \D_{gh}.$$
However, the above property is a consequence of the next two equalities
which can be easily verified by Lemma~\ref{u,v}:

\begin{align*} & (x u_g)y = x (u_g y) \;\;\;\;\; \forall
x\in \B \B_{g\m}, y \in \B_{g\m} \B,\\
& (x v_g)y = x (v_g y) \;\;\;\;\; \forall x\in \B \B_{g} , y \in  \B_{g} \B.\\
\end{align*} Thus $$w_{g,h} \in \M(\D_g \D_{gh}).$$

Finally, because all involved maps are invertible, $w_{g,h}$ is
also invertible with $w\m_{g,h} = u_{gh} v_h v_g.$ \fim \\

\end{section}

\begin{section}{The criteria}\label{sec:criteria}

Now we are ready to formulate the next result.

\begin{teo}\label{criteria} A homogeneously non-degenerate graded algebra
$\B = \bigoplus_{g \in G} \B_g$  is isomorphic as a graded algebra
to a crossed product by a twisted partial action of $G$ if and
only if the following
two conditions are satisfied:\\

(i) $\B_g  \cdot \B_{g\m}  \cdot \B_g    = \B_g  \;\; \forall g \in G;$\\

(ii) for any $g \in G$ there exist $u_g, v_g \in \M (\C_g)$ such that
$u_g v_g = e_{11}$ and  $v_g u_g = e_{22}.$\\

Moreover, if a homogeneously non-degenerate graded algebra $\B$
satisfies (i) and (ii), then with the notation introduced in the
previous sections

$$\Theta = (\{\D_g \}_ {g \in G}, \{{\theta}_g \}_ {g \in G},
 \{w_{g,h} \}_ {(g,h) \in G\times G})$$ form a twisted partial action of
$G$ on $\A = \B_1$ such that\\
$$ \varphi : \bigoplus_{g \in G} \B_g \ni \sum a_g \mt \sum a_gv_g \de _g \in \A \stimes_{\Theta} G$$
is an isomorphism of graded algebras, 
whose inverse is 
$$ {\varphi }\m : \A \stimes_{\Theta} G  \ni \sum a_g \de _g \mt \sum a_gu_g  \in  \bigoplus_{g \in G} \B_g. $$

\end{teo}
\p The ``only if'' part is proved already in
Section~\ref{sec:multipliercross} so suppose that $\B$ is
homogeneously non-degenerate satisfying (i) and (ii). Evidently,
$\D_g = \B_g \B_{g\m}$ $(g \in G)$ are idempotent two-sided ideals
in $\A$ which commute by (ii) of Lemma~\ref{easy}. Thus we have
(i) of Definition~\ref{def1}. We know already by (\ref{isos}) that
$${\0}_{g} : \D_{g\m} \ni x \mt u_gxv_g \in \D_g$$ are isomorphisms of
$k$-algebras, and by Proposition~\ref{w}, $w_{g,h}$ are invertible
multipliers of $\D_g \D_{gh}.$

Obviously, $\D_1 = (\A)^2 = \A,$ by (i). Then evidently $\C_1$ is
the algebra of $2\times 2$-matrices over $\A$ and one can define
the multipliers $e_{12} = \begin{pmatrix} 0 & 1 \\ 0&
0\end{pmatrix} $ and $e_{21} = \begin{pmatrix} 0 &0 \\ 1&
0\end{pmatrix}$ by setting\\

$$\begin{pmatrix} x & y\\ x'  & y' \end{pmatrix} \cdot
\begin{pmatrix} 0 & 1 \\ 0& 0\end{pmatrix} =
\begin{pmatrix} 0 & x\\ 0 &x' \end{pmatrix}, \;\;\;
\begin{pmatrix} 0 & 1 \\ 0& 0\end{pmatrix} \cdot
\begin{pmatrix} x & y\\ x' & y' \end{pmatrix}  =
\begin{pmatrix} x' & y'\\ 0 &0 \end{pmatrix},$$

$$\begin{pmatrix} x & y\\ x' & y' \end{pmatrix} \cdot
\begin{pmatrix} 0 &0 \\ 1& 0\end{pmatrix} =
\begin{pmatrix} y & 0\\ y' &0 \end{pmatrix}, \;\;\;
\begin{pmatrix} 0 & 0 \\ 1& 0 \end{pmatrix} \cdot
\begin{pmatrix} x & y\\ x' & y' \end{pmatrix}  =
\begin{pmatrix} 0 & 0\\ x & y \end{pmatrix},$$\\ where $x,y, x',y' \in \A .$
The multipliers $u_e$ and $v_e$ can be chosen to be equal to
$e_{12}$ and $e_{21}$ respectively. Then $L_{u_1}, L_{v_1} : \A
\to \A, R_{u_1}, R_{v_1} : \A \to \A, \tilde L_{u_1}, \tilde
L_{v_1} :\A \B \to \A \B$ and $\tilde R_{u_1}, \tilde R_{v_1} : \B
\A \to \B \A$ are all identity maps. Consequently, $\0_1:\A \to
\A$ is the identity isomorphism, and $w_{1,g} = u_g v_g = w_{g,1}$
is the identity multiplier of $\D_g^2 =\D_g.$
 This gives (ii) and (v) of Definition~\ref{def1}.

Next observe that
\begin{align*} &\0_g (\D_{g\m} \D_h) = (u_g \D_{g\m} \D_h )v_g =
((u_g \D_{g\m}) \D_h)v_g = (\B_g \D_h) v_g = (\B_g \D_{g\m} \D_h)v_g=\\
& (\B_g \D_h \D_{g\m})v_g  = \B_g \D_h (\D_{g\m})v_g = \B_g \D_h \B_{g\m}
= \B_g \D_h  \B_{g\m} \D_g =  \B_g \B_h \B_{h\m} \B_{g\m} \D_g=\\&
\B_{gh} \B_{h\m g\m} \D_g = \D_{gh} \D_g,
\end{align*} and (iii) of Definition~\ref{def1} is also satisfied.

For arbitrary $x \in \D_{h\m} \D_{h\m g\m}$ we have
$$(w_{g,h} \0_{gh} (x))w\m_{g,h} =
(u_g u_{h} (v_{gh} \circ u_{gh} (x v_{gh}))) u_{gh} v_h v_g =
(u_g u_h (x v_{gh}))u_{gh} v_h v_g,$$ as $v_{gh} \circ u_{gh}$ is identity.
By a direct inspection we see that
$$u_h(x v_{gh}) \in \B_h \B_{h\m g\m} = \B_h \B_{h\m g\m}  \D_{gh}=  \B_{g\m}  \D_{gh} = \D_{g\m} \B_{g\m} \D_{gh}
 \subseteq \D_{g\m} \B \D_{gh}.$$

We also have that
\begin{equation}\label{commuting1} (u_g y)u_{t} = u_g (y u_{t}) \;\;\;
\forall y \in \D_{g\m} \B \D_{t}.
\end{equation} Indeed, write $y = \sum z z' z''$ with
$z \in \D_{g\m}, z' \in \B, z'' \in \D_{t}.$ Then
\begin{align*}& (u_g (z z' z''))u_{t} = (u_g(z)z' z'')u_{t} =
u_g(z)z' (z'' u_{t}) = (u_g z) ((z' z'')u_{t})=\\& u_g(z ((z'
z'')u_{t})) = u_g((z z' z'')u_{t}),
\end{align*}
 which implies the  claimed equality.

Applying (\ref{commuting1}) with $t = gh$ it follows that $$(u_g
(u_h (x v_{gh})))u_{gh} = u_g ((u_h (x v_{gh}))u_{gh}).$$ Now it
is readily seen that
$$ x v_{gh} \in \D_{h\m} \B_{h\m g\m} = \D_{h\m} \B_{h\m g\m} \D_{gh} \subseteq
\D_{h\m} \B \D_{gh},$$
and using  (\ref{commuting1}) again  we have
$$(u_h (x v_{gh}))u_{gh} = u_h (x v_{gh} \circ u_{gh}) =
u_h x, $$ as $v_{gh} \circ u_{gh}$ is identity. Thus we obtain that
$$(w_{g,h} \0_{gh} (x))w\m_{g,h} = (u_g u_h x) v_h v_g.$$ Observe next that
$$u_h x \in  {\B_h} \D _{h\m g\m }= \B _{g\m} \B _{gh}= \D_{g\m}  \B_{h}= \D_{g\m} \B_{h} \D_{h\m}
\subseteq \D_{g\m} \B \D_{h\m}$$
and one can easily check that
\begin{equation}\label{commuting2}(u_g y)  v_h =
u_g (y v_h) \;\;\; \forall y \in \D_{g\m} \B \D_{h\m}.
\end{equation} This yields that
$$(w_{g,h} \0_{gh} (x))w\m_{g,h} = u_g (u_h x v_h) v_g =
\0_g \circ \0_h (x),$$
proving (iv) of Definition~\ref{def1}.

In order to prove that $\Theta$ is a twisted partial action it remains to
check the co-cycle equality (vi) of Definition~\ref{def1}. For any
$x \in  \D_{g\m}  \D_h  \D_{ht}$ we have
\begin{align*} & \0_g (x w_{h,t}) w_{g,ht} =
(u_g (x u_h u_t v_{ht}))v_g u_g u_{ht} v_{ght} = (u_g (x u_h u_t
v_{ht})) u_{ht} v_{ght},
\end{align*} as $v_g u_g$ is identity. We see that

$$x u_h \in \D_{g\m} \D_{ht} \B_h = \D_{g\m} \B_{ht} \B_{t\m},$$
and $$x u_h u_t \in   \D_{g\m} \B_{ht} \D_{t\m} = \D_{g\m} \B_{ht}
\D_{t\m} \D_{t\m h\m} \subseteq \D_{g\m} \B \D_{t\m h\m}.$$ Then
by (\ref{commuting1}) and (\ref{commuting2})
$$u_g(x u_h u_t v_{ht}) = (u_g(x u_h u_t))v_{ht} = ((u_g(x u_h) )u_t )v_{ht},$$ and thus

$$\0_g (x w_{h,t}) w_{g,ht} =
(u_g (x u_h ))u_t v_{ht} u_{ht} v_{ght}= (u_g (x u_h) ) u_t
v_{ght}.$$  Furthermore, since clearly $x \in \D_{g\m}\B \D_h ,$
applying once more (\ref{commuting1}) we obtain that
$$\0_g (x w_{h,t}) w_{g,ht} =  (u_g x) u_h u_t  v_{ght}.$$

On  the other hand,

\begin{align*}  & {\0}_g(x) w_{g,h} w_{gh,t} =
(u_g x) v_g u_g u_h v_{gh}  u_{gh} u_t v_{ght} =
(u_g x)  u_h  u_t v_{ght},
\end{align*} completing the proof of the co-cycle equality.

It remains to show that $\varphi$ is an isomorphism of algebras.
Clearly $a_g v_g \in \D_g$ for $a_g \in \B_g$ so that $a_g v_g
\de_g$ is an element of $\A \stimes_{\Theta} G.$ It is also
evident that $\varphi$ is $k$-linear, as $v_g$ is so. Moreover,
because each $v_g : \B _g \to \D _g $ is bijective with inverse
$u_g,$ it follows that $\varphi $ is bijective with $ \varphi
^{-1} : \sum a_g \de _g  \mt \sum a_g u_g .$ Thus it remains to
check  that $\varphi$ preserves multiplication. Taking $x \in
\B_g$ and $y \in \B_h$ we have
\begin{align*} & \f (x) \f (y) = (xv_g \de_g)\cdot (y v_h \de_h)=
 \0_g({\0}\m_g(x v_g)(yv_h))w_{g,h} \de_{gh}.
\end{align*}
 By Remark~\ref{convenient},
\begin{align*}
&\0_g({\0}\m_g(x v_g)(yv_h))w_{g,h} =
(u_g [ {\0}\m_g (x v_g)y v_h]v_g) u_g u_h v_{gh}=\\&
 (u_g [ {\0}\m_g (x v_g) y v_h]) u_h v_{gh} =
(u_g [ {\0}\m_g (x v_g)]( y v_h)) u_h v_{gh}=\\&
 (u_g [ {\0}\m_g (x v_g)]( y v_h u_h)) v_{gh}=
(u_g [ {\0}\m_g (x v_g)] y ) v_{gh},
\end{align*} as $v_g u_g$ and $v_h u_h$ are  identity.  Furthermore,
\begin{align*}
& (u_g [ {\0}\m_g (x v_g)] y ) v_{gh} = (u_g [ v_g (x v_g u_g)] y
) v_{gh}= (u_g [ v_g x ] y ) v_{gh}=
 (x y)   v_{gh}.
\end{align*} On the other hand,
$$\f (xy) = (x y)v_{gh} \de_{gh}$$ which completes the proof
of the equality $\f(xy) = \f(x) \f(y).$ \fim \\

The  above criteria can be modified under the assumption that each
$D _g$ is an  $s$-unital ring. We remind  that a ring $\R $ is
called {\it left $s$-unital} if for any $r\in \R $ there exists an
element $y \in \R $ with $yr=r.$ Equivalently, for any finite subset $F$ of $ \R $ there exists $y \in \R $ such that $yr =r$ for all $r \in F$ (see \cite[Lemma 2.4]{DdRS}). A left $s$-unital ring is clearly
idempotent. A ring $\R $ is said to be {\it $s$-unital} if it is
both left $s$-unital and  right $s$-unital. We proceed with the
next:

\begin{prop}\label{u,vexistiff}
Let $(\R, {\R}', M, M', \tau, {\tau}')$ be a Morita context with
$s$-unital algebras  $\R $ and ${\R }',$ and surjective $\tau $ and
${\tau }',$ such that 
$ M' \R = {\R }'M' = M'$ and $ M {\R }' = \R M = M.$  Then the following are equivalent:\\

 \noindent (i) There exist $u, v$ in $\M (\C )$ such that  $uv =e_{11}$ and $vu
 =e_{22}.$\\

 \noindent (ii) There are maps 
$$\psi : \R \to M, \; \; \; \mbox{and} \; \; \;  {\psi }' : \R ' \to M, $$
 such that $\psi$ is an isomorphism of left $\R $-modules, ${\psi }'$ is an isomorphism of right ${\R }'$-modules, and
 $$(r {\psi}) r' = r ({\psi }' r') $$ for all $r\in \R $ and $r' \in \R '.$
\end{prop}
\p The algebra $\C $ is easily seen to be idempotent, so (i)
$\Rightarrow $ (ii)  follows by Lemma~\ref{u,v}.

For the converse, observe first that one has the isomorphism $M'
\otimes _{\R } \R \cong M'$ given by $m'\otimes  r \mt m'r.$
Indeed, the map is clearly surjective in view of $M' \R = M'.$
Suppose that $\sum m'_i \otimes r_i \mt 0,$ i. e. $\sum m'_ir_i
=0.$ Then  there exists $r \in \R $ with $ r_i
r = r_i$ for all $i,$ and $ \sum m'_i \otimes r_i = \sum m'_i
\otimes r_i r = \sum m'_i r_i \otimes r = 0,$ as claimed. One
similarly has the isomorphism ${\R }' \otimes _{{\R }'} M' \cong
M'$ defined by $r' \otimes m' \mt r'm'.$ Moreover, the maps $M
\otimes _{{\R }'} M' \to \R ,$ and $M' \otimes _{\R } M \to {\R
}'$ given by the context are also isomorphisms. For the maps
$\tau, \tau '$  are surjective by the assumption. Suppose that
$\tau (\sum m_i \otimes m'_i ) =0.$ In view of $M' \R  = M'$  there exists $r \in \R$ with $ m'_i r
= m'_i $ for all $i.$ Write $r = \tau (\sum _j x_j \otimes x'_j) $
with $x_j \in M, x'_j \in M'.$ Then
\begin{align*} & \sum _i m_i
\otimes m'_i = \sum _i m_i \otimes m'_i r = \sum _i \sum _j  m_i \otimes m'_i \tau (x_j \otimes x'_j ) = \\
& \sum _i \sum _j  m_i \otimes {\tau }' (m'_i \otimes x_j)  x'_j =
\sum _i \sum _j  m_i  {\tau }' (m'_i \otimes x_j) \otimes  x'_j
=\\ & \sum _j \sum _i \tau ( m_i  \otimes m'_i)  x_j \otimes x'_j
= 0,
\end{align*} which shows that $\tau $ is injective, and thus it is an isomorphism.  The verification for  ${\tau }'$
is similar.

 We shall use notation consistent with that of Lemma~\ref{u,v}, so    write
$$ L_u: {\R }' \ni r' \mapsto {\psi }'r' \in M,
 \; \; \;  \mbox{ and } \; \; \;
R_u: \R \ni r \mapsto r{\psi }   \in M .  $$ The composition
\begin{equation}\label{compos1} M' = {\R }' \otimes _{{\R }'} M' \stackrel{ {\psi }'\otimes 1}{\longrightarrow } M
\otimes _{{\R }'} M' \to \R
\end{equation} 
clearly possesses an inverse and thus it is an isomorphism of right $\R $-modules which we  denote  by $L_u$
(so that that as in Lemma~\ref{u,v} $L_u$ stands for the pair $ M'
\to \R , {\R }' \to M $). Use $R_u$ to denote the following
composed map:
\begin{equation}\label{compos2}
   M' = M'  \otimes _{\R } \R  \stackrel{ 1\otimes \psi
}{\longrightarrow } M' \otimes _{\R } M \to {\R }'
 \end{equation} which is an isomorphism of left ${\R }'$-modules
(thus $R_u$ stands now for the pair $\R \to M, M' \to {\R }'$).
Denote the inverses of the above maps by $L_v: \R \to M', \; M\to
{\R }' $ and $R_v : M\to \R , \; {\R }' \to M'.$ Next define the
multiplier $u$ of $\C $ by setting

\begin{equation*} u \begin{pmatrix}r&m\\m'&r'\end{pmatrix} =
\begin{pmatrix} um'& ur'\\0&0\end{pmatrix} \; \;
\mbox{ and }  \; \;  \begin{pmatrix}r&m\\m'&r'\end{pmatrix} u =
\begin{pmatrix} 0& r u\\0&m' u\end{pmatrix},
\end{equation*} where $um' := L_u(m'), \;  ur' := L_u(r'), \;  ru := (r)R_u , \;  m' u :=
(m')R_u.$ One needs to make sure that this in fact  defines a
multiplier. By Remark~\ref{uv(iff)} one has  to check the
equalities in the three columns involving $u$ in Lemma~\ref{u,v}.
The fact that our four maps $L_u, R_u$ are module morphisms give
four of the equalities, i. e. $u (m' r) = (um' )r, u(r'_1 r'_2) =
(u r'_1) r'_2, (r_1 r_2)u = r_1 (r_2 u), (r' m' ) u = r' (m' u),$
whereas the property $( r \psi )r' = r ({\psi }' r')$ means $(r u)
r' = r(u r').$ Moreover, (\ref{compos1}) and (\ref{compos2}) give
 $u(r' m') = (ur')m'$ and $(m'r)u = m'(ru).$  It is a straightforward
exercise to verify the remaining five equalities, and bellow  we
check one of them: $u(m'm) = (um')m.$ Write $m' = \sum r'_i m'_i$
in view of $M' = {\R }' M'.$ Then using the equalities $u(r' m') =
(ur')m'$ and $u(r'_1 r'_2) = (ur'_1)r'_2)$ pointed out above, we
have  $(u m')m = (u \sum r'_i m'_i )m= \sum ((ur'_i)m'_i)m = \sum
(ur'_i)(m'_i m)= \sum u(r'_i (m'_i m)) = \sum u((r'im'_i)m) =
u(m'm),$ as desired.

Analogously, the equalities
\begin{equation*} v \begin{pmatrix}r&m\\m'&r'\end{pmatrix} =
\begin{pmatrix} 0& 0\\v r&v m \end{pmatrix} \; \;
\mbox{ and }  \; \;  \begin{pmatrix}r&m\\m'&r'\end{pmatrix} v =
\begin{pmatrix} m v& 0\\r' v& 0 \end{pmatrix},
\end{equation*} where $vr  := L_v(r), \;  vm  := L_v(m), \;  m v := (m)R_v , \;  r' v :=
(r')R_v,$ define a multiplier $v$ of $C ,$ and it is readily seen
that $uv = e_{11}, vu = e_{22}.$ \fim \\


\begin{remark}\label{remove} In view of Remark~\ref{assumptions} one directly
verifies that the above proposition remains valid if we remove from the assumptions one of the pairs of equalities $\R M = M {\R }' =M$ or $ {\R }' M' =  M'{\R } =  M'.$
\end{remark}

A left module $M$ over a (non-necessarily unital) ring $\R $ is
called {\it torsion-free} (or non-degenerate) if $$x\in M, \; \;
\R x = 0 \Longrightarrow x=0.$$

On easily obtains  the following:

\begin{lema}\label{s-unital}
Let $\R$ be a left $s$-unital ring and let $M$ be a left
$\R$-module such that $M=\R M$. If $x\in M$ then $x\in \R x.$ In
particular, $M$ is torsion-free.
\end{lema}


Now we can state the next:\\

\begin{teo}\label{criteria2} Suppose that $\B = \bigoplus_{g \in G} \B_g$ is a $G$-graded algebra
such that $\D _g$ is an $s$-unital ring for each $g \in G.$ Then
$\B $ is isomorphic as a graded algebra to a crossed product by a
twisted partial action of $G$ on $\B _1$  exactly when the
following
two conditions are satisfied:\\

\noindent (i)  \hspace*{40mm} $\B_g  \cdot \B_{g\m}  \cdot \B_g    = \B_g  \;\; \; \; \; \forall g \in G;$\\

\noindent (ii)  for any $g \in G$ there exist maps
$${\psi }_g:  \D _g \to \B _g \; \; \; \mbox{ and }\; \; \; {\psi
}'_g : {\D }_{g\m } \to \B _g,$$ such that $\psi _g$ is an
isomorphism of left $\D _g $-modules, ${\psi }'_g$ is an
isomorphism of right ${\D }_{g\m }$-modules, and
 $$(d {\psi}_g) d' = d ({\psi }'_{g\m } d') \; \; \; \; \;    \forall d\in \D _g , d' \in \D _{g \m} .$$\\
\end{teo}
\p  By (i) we have $$\D _g \B _g = \B _g = \B _g \D _{g\m }, \; \;
\D _{g\m } \B _{g\m } = \B _{g\m } =\B _{g\m }\D _{g }$$ for each
$g\in G.$ By Lemma~\ref{s-unital}, the left $\D _g$-module $\B
_{g\ }$  is torsion-free. Hence if  $x \in {\B}_g$ and ${\B}_{g\m}
x =0$ then ${\D}_g x=0,$ and consequently $x=0.$ Similarly, $x \in
{\B}_{g\m}, x{\B}_g =0$ implies $x=0,$ and thus  $\B $ is
homogeneously non-degenerate. Our result then follows from
Theorem~\ref{criteria} and Proposition~\ref{u,vexistiff}. \fim

\begin{remark} If $\B = \A \stimes _{\Theta}G$ is a crossed
product by a twisted partial action, then $\psi _g$ is the map $\D
_g \ni a \mt a \de _g \in \B_g$ and ${\psi }'_g$ is $\D _{g\m }
\ni a \mt \0 _g (a) \de _g  \in \B _g .$ If $\Theta $ is global
and $\A $ has $1_{\A }$ then each $\D _g = \A $ and identifying
$\A $ with $\A \de _1 \subseteq \B $ we see that  $\psi _g$ is the
multiplication by $1_{\A } \de _g$ from the right whereas ${\psi
}'$ is multiplication by the same element from the left.\\
\end{remark}

\end{section}

\begin{section}{Constructing some finitely determined maps}\label{sec:fdmaps}

Our next aim is to obtain an algebraic version of Theorem~7.3 from
\cite{E0} for which we do some preparatory work in this section.
As it is mentioned in the Introduction, we view the
$C^{\ast}$-property of having a countable ap\-pro\-xi\-ma\-te
identity algebraically as that of possessing a countable set of
local units. We say that an algebra $\R $  has {\it local units}
if for any $r \in \R$ there exists an idempotent $e \in  \R $ such
that $r \in e \R e.$ Equivalently, any finite subset of $\R$ is
contained in a subring of the form $e \R e$ for some $e = e^2 \in
\R $ (see Remark~\ref{localunits} below). The algebra $\R $
possesses a countable set of local units exactly when there exists
an increasing  sequence $e_1 \leq e_2 \leq \ldots $ of idempotents
in $\R $ with $\R = \cup _{i=1}^{\infty } e_i \R e_i .$  Here $e_1
\leq e_2 $ means that $e_1 e_2 =e_2e_1 = e_1 .$ The elements $e'_1
=e_1$ and $e'_i = e_i -e_{i-1}$ $(i \geq 2)$ are pairwise
orthogonal idempotents and $\R = \bigoplus_{i \in \Nn} \R e'_i =
\bigoplus _{i \in \Nn} e'_i \R .$  More generally, we say that
$\R$ possesses {\it orthogonal local units} if there exists a set
of (non-necessarily central) pairwise orthogonal idempotents $E$
in $\R$ such that

\begin{equation}\label{ortolocal1}
 \R = \bigoplus_{e\in E}\; \R e = \bigoplus_{e\in E} \; e \R .
\end{equation}

\noindent Rings $\R $ with (\ref{ortolocal1}) are also called
rings with enough idempotents (see \cite{Fuller}). As it is shown
by P. N. \'{A}nh and L. M\'{a}rki \cite{AnhMarki}, every Morita
equivalence class of rings with local units contains rings with
orthogonal local units.

\begin{remark}\label{localunits} We observe that if $\R $ is a
ring such that for any $r\in \R$ there exists an idempotent $e\in
\R$ with $r \in e \R e,$ then for any finite subset $\{ r_1, r_2,
\ldots r_m \}  \subseteq \R $ there exists an idempotent $e\in \R$
such that $r_i \in e \R e$ for all $i=1, \ldots m.$ Indeed, our
ring $\R $ is  right $s$-unital and thus there exists $y \in \R$
such that $r_iy=r_i$ for all $i=1,\ldots ,m.$ Take an idempotent
$f\in \R$ with $yf=y.$ Then $r_i = r_iyf = r_if,$ $i=1,\ldots m.$
Analogously, since $\R$ is left $s$-unital, there exists $x\in \R$
with $xr_1 =r_1, \ldots, xr_m =r_m$ and $xf=f.$ Take an idempotent
$e   \in \R$ such that $ex=x.$  Then $r_i = er_i,$ $i=1,\ldots m$
and $f=ef.$ Setting now $e' = e+f-fe$ it is directly verified that
$e'$ is an idempotent with $r_1,\ldots r_m \in e'\R e'.$
\end{remark}

We shall need some information on the Morita Theory for idempotent
rings. The reader is referred to \cite{GS} for the details. Given
an idempotent ring $\R ,$ a left $\R$-module is said to be  {\it
unital} if $\R M = M.$  The category of left (respectively, right)
unital torsion-free modules is denoted by $\R$-${\rm mod}$
(respectively, ${\rm mod}$-$\R $). Two idempotent rings $\R$ and
${\R}'$ are called Morita equivalent if the categories $\R$-${\rm
mod}$ and ${\R}'$-${\rm mod}$ are equivalent. This happens exactly
when there is a Morita context $(\R, {\R}', M, M', \tau, {\tau}')$
with unital $_{\R} M, M _{{\R}'}, _{{\R}'} M', M'_{\R }$ and
surjective $\tau : M\otimes_{\R'} M' \to \R,$  $\tau' :
M'\otimes_{\R} M \to \R'.$
Then the categories ${\rm mod}$-$\R $ and ${\rm mod}$-${\R}'$ are also equivalent and vice versa.\\



When dealing with modules over rings without identity, the concept
of a finitely ge\-ne\-ra\-ted module is refined by using the
categorical definition with respect to a family of generators (see
\cite[p. 72]{M} ). For rings with orthogonal local units, however,
we will be satisfied with the ``usual'' definition, so that we say
that a left $\R $-module $N$ is {\it finitely generated } if there
exist finitely many $y_i \in N$ such that $N = \sum _i \R y_i.$
Note that $\R $ in general does not need to be finitely generated as
a module over itself. Note also that a finitely generated $\R
$-module is necessarily unital, as $\R $ is idempotent.\\

We shall use the next  easy fact  which also explains  that our
definition of finitely generated modules is equivalent to the
categorical one with the most naturally chosen family of
generators.

\begin{lema}\label{f-generated}
Let $\R$ be a ring with a set of orthogonal local units $E$ and
set $\mathcal{E}= \{ {\R} e : e \in E \}.$ A left $\R$-module $M$
is $\mathcal{E}$-finitely generated if and only if there exists a
finite number of elements $y_1,\dots, y_t$ of $M$, such that
$M=\sum \R y_t$.
\end{lema}

\p  For the reader's convenience we remind that according to \cite[p. 72]{M},  a left $\R $-module $M$  is said to be finitely generated with respect
to the generating set $\mathcal{E}$ if there exists a finite
subset $E_0 \subseteq E$ such that $M$ is an epimorphic image of a
finite direct sum of modules each of which is isomorphic to some
${\R}e$ with $ e\in E_0.$  Now our assertion easily follows by observing that for any $x\in
M$, there exists a finite subset $E'\subseteq E$
    such that $x\in \sum_{e\in E'}\R e x$. \fim \\

Let $P$  be a left $\R $-module with given fixed decomposition
$$P = \bigoplus P_i$$ into a (non-necessarily finite) direct sum of finitely generated
submodules. By an {\it admissible summand} of $P$ we shall mean a
direct summand of $P$ which is  a finite sum of some of the
$P_i$'s. The $P_i$'s shall be called {\it admissible components}
of $P.$ Let $Q$ be
 also a left $\R $-module with the fixed decomposition
$$Q = \bigoplus Q_j,$$ with finitely generated $Q_j$'s which are taken to be the
admissible components of $Q.$

\begin{defi}  An $\R $-module  homomorphism $\psi : P \to Q $  shall be
called
 finitely determined if for any admissible  summand  $L_1$ of $Q$   there exists
 an admissible summand $L_2$ of $P$ such that
$$(x) \psi \circ {\lambda}_1 = (x){\lambda}_2 \circ \psi |_{L_2} \circ {\lambda}_1$$
for all $x \in P,$ where ${\lambda}_1 :Q\to L_1$ and ${\lambda}_2 :P \to L_2$ are the projections and $\psi
|_{L_2}$ is the restriction.
\end{defi}

Evidently the composition of finitely determined maps is finitely determined, and when checking that $\psi$ is
finitely determined is enough to take $L=Q_j$ with arbitrary $Q_j.$

For the module $P \bigoplus Q$ the admissible components will be those of $P$ and $Q,$ and the same will be used
for infinite
direct sums of modules with fixed admissible components.\\

In order to construct some finitely determined maps we shall fix the following notation. The rings  $\R$ and ${\R}'$ both will have orthogonal local units, respectively denoted by $E$ and $F.$ Thus $\R
$ satisfies (\ref{ortolocal1}) and similarly for ${\R}'$ we have

\begin{equation}\label{ortolocal2}
{\R}' = \bigoplus_{f\in F}\; {\R}' f = \bigoplus_{f\in F} \; f {\R}'.
\end{equation}

The cardinalities of $E$ and $F$ are arbitrary. Suppose further that $\R$ and ${\R}'$ are Morita equi\-va\-lent.
This means that there exists a Morita context $(\R, {\R}', M, M', \tau, {\tau}')$ with surjective maps (in fact,
by \cite[Proposition 2.4]{GS} the maps $ \tau, {\tau}'$ are bijective).  Moreover, by Lemma~\ref{s-unital}, the modules
$_{\R} M, M _{{\R}'}, _{{\R}'} M', M'_{\R }$ are all torsion-free and by \cite[Theorem 2.7]{GS}, the functors
${{M}'\otimes _{\R}} \; \underline{\;\;} : \R$-${\rm mod} \to {\R}'$-${\rm mod} $ and $M \otimes _{{\R}'} \;
\underline{\;\;} : {\R}'$-${\rm mod} \to {\R}$-${\rm mod} $ are inverse category equivalences, which are unique
up to natural isomorphisms. Note that our category equivalences preserve short exact sequences. Write $x \cdot
x' = \tau (x\otimes x')$ and $x' \cdot x = \tau' (x' \otimes x)$ where $x \in M, x' \in M'.$   Since $M = {\R }M
= M{{\R}'}$  and $M' = {\R }'M'
= M'{{\R}},$  it follows by (\ref{ortolocal1}) and (\ref{ortolocal2}) that
\begin{equation}\label{direct1}
M= \bigoplus _{e \in E}eM \, \, = \, \,  \bigoplus _{f \in F}Mf \, \, \, \, \, \,\mbox{and} \, \, \, \, \, \,
M'= \bigoplus _{f \in F}fM' \, \, = \, \,  \bigoplus _{e \in E}M' e.
\end{equation}

For any  $f \in F$ write
\begin{equation}\label{f}
 f =
\sum_{i=1}^{n_f} x'_i \cdot x_i,
\end{equation} where $x_i =
x^{(f)}_i \in M,$ and $ x'_i = x'^{(f)}_i  \in M'.$  One evidently may suppose that  $x_i = x_i f$, and $f x'_i =x'_i$ for all $i.$\\

  Observe next the following easy facts.

\begin{lema}\label{fingen} With the above notation we have:\\

 (i) For each $e \in E$ and $f \in F$ the modules $_{\R}(Mf),$
and $ _{{\R}'}(M'e)$ are finitely generated.\\

    (ii)  Finitely generated modules are preserved under the equivalences ${{M}'\otimes _{\R}} \;
\underline{\;\;}$ and $M \otimes _{{\R}'} \; \underline{\;\;}.$
\end{lema}

\p (i) For any $y \in M$ by (\ref{f})  one has $yf =\sum_i (yx'_i) x_i .$   Since $yx'_i \in \R ,$  it follows that $Mf = \R x_1 + \ldots + \R x_{n_f},$ so that   the
left $\R$-module $Mf$ is
finitely ge\-ne\-ra\-ted. Similarly, so too is $ _{{\R}'}(M'e).$\\

(ii) Follows easily by using  (i) and Lemma~\ref{f-generated}.
\fim\\

 We proceed by constructing some homomorphisms.  For each $f\in F$, by (\ref{f}) we may define the map

$${\pi}_f:{\R}^{n_f} \ni (r_1, r_2,  \ldots  r_{n_f}) \mapsto
\sum r_i x_i \in Mf,$$  where $x_i=x^{(f)}_i$ and $x'_i=x'^{(f)}_i.$ This is an epimorphism of left
${\R}$-modules, as  $yf = {\pi}_{f} (yx'_1, yx'_2, \ldots , yx'_{n_f})$  for any $y \in M .$ The map
$${\rho}_f : Mf \ni yf \mapsto ( yfx'_1, yfx'_2, \ldots yfx'_{n_f}) \in
{\R}^{n_f},$$ obviously splits ${\pi}_f,$ i.e ${\rho}_f \circ
{\pi}_f = id .$  Then we have the  exact sequence of left $\R
$-modules and homomorphisms
\begin{equation*}\label{exact1}
0 \longrightarrow  K_{f} \stackrel{{\mu}_f}\longrightarrow
{\R}^{n_f} \stackrel{{\pi}_f}\longrightarrow  Mf \longrightarrow
0,
\end{equation*} in which $K_f$ is the kernel of ${\pi}_f$ and ${\mu
}_f$ is its imbedding, and also the exact sequence of the
splitting maps
\begin{equation*}\label{exact2}
0 \longleftarrow K_{f} \stackrel{{\tau}_f}\longleftarrow
{\R}^{n_f} \stackrel{{\rho}_f}\longleftarrow  Mf \longleftarrow 0
\;\; ({\mu}_f \circ {\tau}_f =id).
\end{equation*}\\  Note that since $x'_i = fx'_i,$ one can
consider ${\rho}_f$ to be defined on $M$ and $x {\rho }_f =
(xf){\rho }_f$ for all $x \in M.$

Let $X$ be an infinite set of indices whose cardinality is bigger
or equal to those of $E$ and $F.$ Completing $E$ or $F$ with
zeros, we may assume that $|X| = |E| = |F|.$  Let further
${\R}^{(X)}$ be the direct sum of copies of ${\R}$ indexed by the
elements of $X$ and write ${\R}^{(X)} = \bigoplus _{f \in F}
{\R}^{n_f}.$ Here we assume ${\R}^{n_f}= 0$ if $f =0.$ Define the
map
$$\pi : \bigoplus _{f \in F} {\R}^{n_f} \ni ({\bar r}_f)  \mapsto
\sum _{f\in F}({\bar r}_f){\pi}_f \in M,$$ here ${\bar r}_f \in
{\R}^{n_f}.$ By  (\ref{direct1})  $\pi $ is an epimorphism of left
${\R}$-modules. Since $x= \sum _{f\in F} xf,$ we see that  $\pi $
is split by the homomorphism
$$\rho : M \ni x \mapsto (x{\rho}_f) \in \bigoplus _{f\in F} {\R}^{n_f},$$ i.e.
$\rho \circ \pi = id.$
Similarly,  the maps  ${\mu}_f,$  ${\tau}_f$ $(f \in F)$ determine the left $\R $-module
homomorphisms
$\mu : K \to  {\R}^{(X)}$ and  $\tau : {\R}^{(X)} \to K,$ where $K = \bigoplus _{f \in F} K_f,$
resulting in the exact
sequences
\begin{equation}\label{exact3}
0 \longrightarrow  K \stackrel{{\mu}}\longrightarrow
{\R}^{(X)} \stackrel{{\pi}}\longrightarrow  M \longrightarrow 0 \; \; \mbox{and}\; \;
0 \longleftarrow K \stackrel{{\tau}}\longleftarrow
{\R}^{(X)} \stackrel{{\rho}}\longleftarrow  M \longleftarrow 0,
\end{equation} where ${\mu} \circ {\tau} =id.$\\

Taking ${\R}'$ instead of $\R $ and $M'$ instead of $M$ in
(\ref{exact3}), we obtain that  there are exact sequences of left
${\R}'$-modules and  homomorphisms

\begin{equation}\label{exact3.5}
0 \longrightarrow K' \stackrel{{\mu}'}\longrightarrow
{{\R}'}^{(X)} \stackrel{{\pi}'}\longrightarrow  M' \longrightarrow
0 \; \; \mbox{and}\; \; 0 \longleftarrow K'
\stackrel{{\tau}'}\longleftarrow {{\R}'}^{(X)}
\stackrel{{\rho}'}\longleftarrow M' \longleftarrow 0,
\end{equation} with ${\rho}' \circ {\pi}' =id, {\mu}' \circ {\tau}'
=id.$ Applying the equivalence functor ${M\otimes _{{\R}'}} \;
\underline{\;\;}$ we come to the exact sequences of left
${\R}$-modules and homomorphisms:

\begin{equation}\label{exact4}0 \longrightarrow  \tilde{K} \stackrel{\tilde{\mu}}\longrightarrow
{M}^{(X)}
\stackrel{\tilde{\pi}}\longrightarrow  {\R} \longrightarrow 0 \;
\; \mbox{and}\; \; 0 \longleftarrow \tilde{K}
\stackrel{\tilde{\tau}}\longleftarrow M^{(X)}
\stackrel{\tilde{\rho}}\longleftarrow  {\R} \longleftarrow 0,
\end{equation} with $\tilde{\rho} \circ \tilde{\pi} =id, \tilde{\mu} \circ
\tilde{\tau} =id,$ where $ \tilde{K}$ is the submodule of
${M}^{(X)}$ obtained from ${M\otimes _{{\R}'}} K'$ by using
${M\otimes _{{\R}'}} {\R}' \cong _{\R} M,$ and $\tilde{\mu}$ is
the embedding
of $ \tilde{K}.$\\

We shall fix for the rest of this section the decomposition for the modules $_{{\R}} M,$ $ _{{\R}'} M',$ $ _{\R}
{\R},$ $ _{{\R}'} {\R}',$ $ _{{\R}} K,$ $ _{{\R}'} K',$ and $ _{{\R}} \tilde{K}.$ For $_{\R} {\R}$ and
$_{{\R}'} {\R}'$ take as fixed the decompositions

$$\R = \bigoplus_{e\in E}\; \R e \, \, \, \, \,
\mbox{and} \, \, \, \, \,  {\R}' = \bigoplus_{f\in F}\; {\R}' f,$$

\noindent so that the admissible components of $\R $ are the $\R
e$'s and those of ${\R}'$ are the  ${\R}' f$'s.\\

We fix now the decompositions for  $_{{\R}} M$ and  $_{{\R}'} M'$ :

$$   M=  \bigoplus _{f \in F}Mf  \, \, \, \, \,
\mbox{and} \, \, \, \, \,    M'=  \bigoplus _{e \in E}M' e.$$\\

As to $_{{\R}} K,$  since $M=  \bigoplus _{e \in E}e M,$ there exists an idempotent $e^{(f)} \in \R$ such that
$e^{(f)} = e^{(f)} x^{(f)}_j$ for all $j =1, \ldots n_{f},$ which  is a (finite) sum of some $e$'s from $E.$
Clearly $$K^{(f)} = \bigoplus _{e\in E, e e^{(f)}=0} ({\R }e)^{n_{f}} \subseteq K_{f},$$ and we have $K_{f} =
K^{(f)} \bigoplus L^{(f)} $ with $L^{(f)} = K_{f} \cap  ({\R }e^{(f)})^{n_{f}}.$ Since $K_{f}$ is a direct
summand of ${\R}^{n_{f}},$ the left $\R $-module $L^{(f)} $ is a direct summand of $({\R }e^{(f)})^{n_{f}}.$
Because  $({\R }e^{(f)})^{n_{f}}$ is finitely generated, so too is $L^{(f)} .$ Then the sum

\begin{equation}\label{admisK}
 K = \bigoplus _{f \in F} \; \; \, \, ( L^{(f)} \bigoplus  _{e\in
E, e e^{(f)}=0} ({\R }e)^{n_{f}} )
\end{equation}

\noindent is a decomposition of $K$ into finitely generated direct summands, and we fix the mo\-du\-les
$L^{(f)}$ and ${\R }e$ with $f \in F$ and $e e^{(f)} =0, e\in E,$ as the admissible components of $K.$

Working similarly with ${\R}',$  $M'$ and $K',$ we define for each
$e\in E$ the idempotent $f^{(e)}$ and the module $L'^{(e)},$ so
that

\begin{equation}\label{admisK'}
 K' = \bigoplus _{e \in E} \; \; \, \, ( L'^{(e)} \bigoplus
_{f\in F, f f^{(e)}=0} ({\R }'f)^{n_{e}} )
\end{equation}

\noindent is a decomposition of $K' $ as a direct sum of finitely generated submodules, and the modules $L'^{(e)}$
and ${\R }'f$ with $e \in E$ and $f f^{(e)} =0, f\in F,$ are defined to be the admissible components of $K'.$

Applying the  functor ${M\otimes _{{\R}'}} \; \underline{\;\;}$ to
the decomposition of $K'$ we obtain the decomposition:

\begin{equation}\label{admistildeK}
  \tilde{K} = \bigoplus _{e \in E} \; \; \, \, (
{\tilde{L}}^{(e)} \bigoplus _{f\in F, f f^{(e)}=0} (Mf)^{n_{e}} )
\end{equation}

\noindent of $ \tilde{K}$ into a direct sum of submodules, where ${\tilde{L}}^{(e)}$ is the submodule of $M^{(X)}$
obtained from ${M\otimes _{{\R}'}}   L'^{(e)}$ by using ${M\otimes _{{\R}'}} {{\R}'} \cong M.$ We know already
that each $M f$ is finitely generated and, moreover,  by (ii) of Lemma~\ref{fingen}, so too is each
${\tilde{L}}^{(e)}.$ Thus the modules $Mf$ and ${\tilde{L}}^{(e)}$ with $e \in E$ and $f f^{(e)}=0, f \in F,$
can be defined as the admissible components of $\tilde{K},$ so that the functor ${M\otimes _{{\R}'}} \;
\underline{\;\;}$ preserves the admissible components when passing  from (\ref{admisK'})
to (\ref{admistildeK}).\\

Now we are ready to state the next:

\begin{lema}\label{findet} All homomorphisms in  (\ref{exact3}) and  (\ref{exact4}) are finitely
determined.
\end{lema}
\p Taking $L_1 =  Mf$ with arbitrary  $f \in F$ and  $L_2 =
({\R}e^{(f)})^{n_{f}}, $ it is readily seen that  $\pi $ is
finitely determined. Next, any admissible summand  $L_1$  of
${\R}^{(X)}$ is contained in some  ${\R}^{n_{f}},$ so that taking
 $L_2 =  Mf$ we immediately have
that  $\rho $ is finitely determined. On the other hand, it follows from (\ref{admisK}) that each admissible
summand of $K$ is either equal to an admissible one of ${\R}^{(X)}$ or  contained in a finite direct sum of
them. It follows that both $\mu $ and $\tau $  are finitely determined.

Analogously all maps in  (\ref{exact3.5}) are finitely determined,
and since the admissible summands of (\ref{admistildeK}) are
defined so that they are  preserved by  ${M\otimes _{{\R}'}} \;
\underline{\;\;} \,,$ all maps in  (\ref{exact4}) are also
finitely determined. \fim \\

We apply the above lemma to give the next:

\begin{prop}\label{Eilenberg} There is a finitely determined
isomorphism of left $\R $-modules\\ $\psi :  {{\R}}^{(X)} \to
M^{(X)}$ whose inverse is also finitely determined.
\end{prop}
\p Observe first that by Lemma~\ref{findet} the following are
finitely determined isomorphisms of left $\R $-modules:
$$\pi \oplus \tau : {{\R}}^{(X)} \la M \oplus K, \;\;
\tilde{\pi} \oplus \tilde{\tau} : M^{(X)} \la \R \oplus
\tilde{K},$$ whose inverses are also finitely determined. Then
using these maps and rearrangements of indices we apply  the
Eilenberg's trick:
\begin{align*} &
{\R}^{(X)} \cong ({\R}^{(X)})^{(X)} \cong (M \oplus K)^{(X)} \cong
M^{(X)} \oplus K^{(X)} \cong M^{(X)} \oplus M^{(X)} \oplus
K^{(X)}\\& M^{(X)} \oplus (M \oplus K)^{(X)} \cong M^{(X)} \oplus
({\R}^{(X)})^{(X)} \cong  M^{(X)} \oplus {\R}^{(X)},
\end{align*} and similarly
 $$M^{(X)} \cong ({M}^{(X)})^{(X)} \cong (\R \oplus \tilde{K})^{(X)}
 \cong \ldots \cong {\R}^{(X)} \oplus {M}^{(X)}. $$ It is directly seen
 that each step is a finitely determined
isomorphism and thus  so too is their composition $\psi :
{\R}^{(X)} \to M^{(X)},$ as well as
${\psi}^{-1}.$ \fim \\

  It seems to be  interesting  to give here  some easy examples of
  non-finitely determined isomorphisms.

\begin{ex}\label{nonfindet}  Let  $\phi :  {\Rr}^{(\Nn)} \to {\Rr}^{(\Nn)} $
 be defined  by $e_1 \mapsto e_1, e_i \mapsto e_1 + e_i, i \geq 2,$
where the $e_i$'s are the canonical basis of ${\Rr}^{(\Nn)}.$ Then
$\phi  $ is an ${\Rr}$-isomorphism whose inverse is determined by  $e_1
\mapsto e_1, e_i \mapsto  e_i - e_1, i \geq 2.$ Taking the
canonical copies of $\Rr $ as the admissible components, it is
immediately seen that $\phi $ is not finitely determined, as well
as its inverse.
\end{ex}

\begin{ex}\label{nonfindet} Using notation from the above example define   $\phi :  {\Rr}^{(\Nn)} \to {\Rr}^{(\Nn)} $
 by $e_1 \mapsto e_1, e_i \mapsto e_{i-1} + e_i, i \geq 2.$  Then $\phi$ is an ${\Rr}$-isomorphism whose inverse is given by 
$e_i \mapsto 
\sum_{j=0}^{i-1} (-1)^{j}e_{i-j}.$  It is readily  seen that $\phi $ is
finitely determined, while $\phi^{-1}$ is not.
\end{ex}

\end{section}

\begin{section}{The stable characterization of graded algebras}\label{sec:stab}

 Recall that given an algebra $\R$ and a (non-necessarily finite) set of
indices $X,$ we denote by ${\rm RFMat}_{X}(\R)$ the algebra of all $ X  \times X$-matrices over $\R$ which are
row-finite and by ${\rm FMat}_{X}(\R)$ the algebra of all $ X  \times X$-matrices over $\R$ with only a finite
number of non-zero entries.

\begin{lema}\label{matrices}
Let $\R$ be a ring with a set $E$ of orthogonal local units and $P$ and $Q$ left unital $\R$-modules with
decompositions $P=\oplus _{i\in X}P_i$ and $Q=\oplus_{j\in Y} Q_j$, where $P_i$ and $Q_j$ are finitely
generated. Then \\

(i) ${\rm Hom} _{\R}(P, Q)\cong {\rm RFMat}_{X\times Y}\left({\rm Hom} _{\R}(P_i, Q_j)\right).$\\

(ii)  $ P\cong {\rm FMat}_{(e,i)\in E\times X}(eP_i).$
\end{lema}

\p (i) Since each $P_i$ is a finitely generated module then ${\rm Hom} _{\R}(P_i,\oplus_{j\in Y} Q_j)\cong
\oplus_{j\in Y} {\rm Hom} _{\R}(P_i, Q_j)$.  Thus  ${\rm Hom}
_{\R}(\oplus _{i\in X}P_i,\oplus_{j\in Y} Q_j)\cong \prod_{i\in X}\oplus_{j\in Y} {\rm Hom} _{\R}(P_i, Q_j)\cong
{\rm RFMat}_{X\times Y}\left({\rm Hom} _{\R}(P_i, Q_j)\right).$

(ii) Clearly, ${\rm Hom} _{\R}(\R e, P_i)\cong eP_i$. Then, by (i), ${\rm Hom} _{\R}(\R, P)\cong {\rm
RFMat}_{E\times X}\left(eP_i\right).$ Now $P$ embeds in ${\rm Hom} _{\R}(\R, P)$ as right multiplication; so
that, for each $p\in P$, with $p=\sum p_i$, we have  $p\mapsto \left(ep_i\right)\in {\rm FMat}_{E\times
X}(eP_i)$. Conversely, one may check that $a\in {\rm FMat}_{E\times X}(eP_i)$, comes from the right
multiplication by $p=\sum a(e,i)$.
\fim\\

The following is the main fact that we need about Morita equivalent rings with orthogonal local units.

\begin{teo}\label{uv-for-Morita} Let $\R$ and ${\R}'$  be rings with
orthogonal local units $E$ and $F$ respectively and $X$ be an
infinite  set of indices whose cardinality is bigger  than or
equal to those of $E$ and $F.$ Let further   $(\R , {\R}', M, M',
\tau, {\tau}')$ be a Morita context with unital (and thus, by
Lemma~\ref{s-unital}, torsion-free) $_{\R} M, M _{{\R}'}, _{{\R}'}
M', M'_{\R }$ and surjective $\tau : M\otimes_{\R'} M' \to \R,$
$\tau' : M'\otimes_{\R} M \to \R' ,$ and let $\C$ be the
corresponding linking algebra. There exist  maps $$\Psi : {\rm
FMat}_{X}(\R ) \to {\rm FMat}_{X}( M ), \; \; \; \mbox{and} \; \;
\; {\Psi }' : {\rm FMat}_{X} (\R ' ) \to {\rm FMat}_{X} (M ), $$
 such that $\Psi$ is an isomorphism of left ${\rm FMat}_{X}( \R )$-modules, ${\Psi }'$ is an isomorphism of right
 ${\rm FMat}_{X}({\R }')$-modules, and
 $$(A {\Psi}) A' = A ({\Psi }' A') $$ for all $A\in {\rm FMat}_{X}(\R ) $ and $A' \in {\rm FMat}_{X}(\R ').$
\end{teo}

\p A homomorphism ${\phi}$ of left $\R $-modules ${{\R}}^{(X)} \to
M^{(X)}$ can be viewed as an $X \times X$-matrix with entries in
${\rm Hom} _{\R}(\R, M).$  By (i) of Lemma~\ref{matrices} ${\phi}$
becomes an element $[{\phi}]$ of the set ${\rm Mat}_{X} ({\rm
RFMat}_{(e,f) \in E \times F} (eMf)),$ i. e. $[{\phi}]$ has block
structure the block-rows and block-columns being indexed by the
elements of $X,$ and each block is a matrix from ${\rm
RFMat}_{(e,f) \in E \times F} (eMf).$

 Fixing $i \in X$ and  $e\in E$ we have the
(``thin'') $(i,e)$-row, which has to be finite, as the sum of its
elements is the  image under ${\phi}$ of $e$ from the $i$-copy of
$\R $ in ${{\R}}^{(X)}.$ Now taking the finitely determined
isomorphism ${\psi}$ from Proposition~\ref{Eilenberg}, we have
that  for any $i \in X$ and $f \in F$ the (``thin'')
$(i,f)$-column of $[{\psi}]$ is also finite.  It follows that the
sum of each block-row as well as that of each block-column of
$[{\psi}]$ makes sense, and one may call such matrices
 as {\it row and column summable} $X \times X$-matrices
over   ${\rm Hom} _{\R}(\R, M)$ and denote the set of them by ${\rm RCSumMat}_{X}({\rm Hom} _{\R}(\R, M)).$ It
is evident that each block of $[{\psi}]$ is a row and column finite matrix, so that denoting the set of all row
and column finite matrices from ${\rm Mat}_{(e,f) \in E \times F} (eMf)$ by ${\rm RCFMat}_{(e,f) \in E \times F}
(eMf),$ one can write that  $$[{\psi}] \in {\rm RCSumMat}_{X}({\rm RCFMat}_{(e,f) \in E \times F} (eMf)).$$ It
may be interesting to observe that the sum of each block-row of $[{\psi}]$ is a row-finite matrix whereas that
of each block-column is column-finite.

Now by (ii) of Lemma~\ref{matrices} we have the following ring isomorphisms

\begin{equation*}\label{R}
\e :\R \cong {\rm FMat}_{e,e'\in E}(e \R e')\, \, \, \, \mbox{and} \, \, \, \, \e' :\R' \cong {\rm
FMat}_{f,f'\in E}(f \R f')
\end{equation*}
and a $k$-linear  isomorphisms

$$\z : M \la {\rm FMat} _{e\in E, f\in F} (eMf) \, \, \, \mbox{and} \, \, \,
\z' : M' \la {\rm FMat} _{f\in F, e\in E} (fM'e) . $$

Moreover, using matrix multiplication, one can check that ${\rm
FMat} _{e\in E, f\in F} (eMf)$ has a  structure of a left ${\rm
FMat}_{e,e'\in E}(e \R e'),$  right ${\rm FMat}_{f,f'\in F}(f {\R
}' f')$-bimodule, and similarly ${\rm FMat} _{f\in F, e\in E}
(fM'e)$ also has its corresponding bimodule structure,  so that
the maps $\z$ and $\z'$ are both semilinear, i. e. $$\z (ry) = \e
(r) \z (y), \z (yr') = \z (y) \e' (r'), \z' (r' y' ) = \e' (r' ) \z' (y'), \z'
(y' r) =\z' (y') \e (r),$$ for any $r\in \R , r'\in {\R}' , y\in M, y'
\in M'.$ Extending $\e , \z , {\e}' , {\z}'$ to finite matrices,
the semi-linearity is obviously preserved. Moreover, the maps
$\tau$ and $\tau '$ from the Morita context yield the equalities
\begin{equation*}
 \e (YY') = \z (Y) {\z}' (Y'),\, \, \, \,  {\e}' (Y'Y) = {\z}' (Y') {\z} (Y).
\end{equation*} for all $Y \in {\rm FMat}_X ( M ), Y' \in {\rm FMat}_X ( M' ).$

An element from ${\rm FMat}_X ({\rm FMat}_{e,e'\in E}(e \R e'))$ can be multiplied from the right by $[{\psi}]$
and
\begin{equation*}\label{rpsi1}
{\rm FMat}_X ({\rm FMat}_{e,e'\in E}(e \R e')) \cdot [{\psi} ] \subseteq {\rm FMat}_X ({\rm FMat}_{ E\times F}(e
M f)).
\end{equation*} Then setting $A \cdot [{\psi} ] = {\z}\m
({\e}(A)[{\psi}])$ with $A\in {\rm FMat}_X ( \R ),$ the map
\begin{equation}\label{psi2}
\Psi :{\rm FMat}_X ( \R ) \ni A \mt A \cdot [{\psi} ] \in {\rm
FMat}_X (M)
\end{equation} is clearly $k$-linear and
\begin{equation}\label{psi3}
(A B) \cdot [{\psi} ] = A (B \cdot [{\psi} ])
\end{equation} for all $A, B \in {\rm FMat}_X ( \R ).$


One also defines the $k$-linear map


$${\Psi }' : {\rm FMat}_X ( {\R}' ) \ni A' \mt   [{\psi} ]\cdot A'  \in {\rm FMat}_X
(M),$$ such that


$$[{\psi} ]\cdot
(A'_1A'_2) =  ([{\psi} ]\cdot A'_1)A'_2$$ and $$ (A \cdot [{\psi}
])A' = A ([{\psi} ]\cdot A')$$

\vspace*{3mm}

\noindent  for  any $A\in {\rm FMat}_X ({\R}), A', A'_1, A'_2 \in
{\rm FMat}_X ({\R}').$\\

Now as it was done for  $\psi$,  the isomorphism ${\psi}\m :
M^{(X)} \to {\R}^{(X)},$ may be viewed as a matrix in ${\rm
Mat}_{X}({\rm RFMat}_{(f,e) \in F \times E} ({\rm Hom} _{\R}(Mf,
\R e))),$ by Lemma~\ref{matrices}. The equivalence functor
$M'\otimes_{\R}\;\underline{\;\;}$ transforms
 ${\rm Hom} _{\R}(Mf, \R e))$ into $ fM'e$; so that
${\rm Hom} _{\R}(M, \R))\cong {\rm RFMat}_{(f,e) \in F \times E}
(fM'e).$ As in the case of  $\psi$, it is directly verified that
$$[{\psi}\m] \in {\rm RCSumMat}_{X}({\rm RCFMat}_{(f,e) \in F
\times E} (fM'e)).$$ Also, as it was done for $[\psi ]$ we define
the map

$${\rm FMat}_X ( M ) \ni Y \mt Y \cdot [{\psi}^{-1} ] \in {\rm FMat}_X
(\R ),$$ which is  ${\Psi }^{-1},$ and the map

$$ {\rm FMat}_X ( M ) \ni Y \mt   [{\psi}^{-1} ]\cdot Y  \in {\rm FMat}_X
( {\R }'),$$ which is the inverse of ${\Psi }'.$  \fim \\

\begin{remark}\label{remark:u,v} Under the assumptions in the above result,   Proposition~\ref{u,vexistiff} guarantees that  there exist
multipliers $u,v$ of $${\rm FMat}_{X} (\C) \cong
\begin{pmatrix} {\rm FMat}_{X} (\R)   &{\rm FMat}_{X} (M)\\ {\rm FMat}_{X} (M')   &{\rm FMat}_{X} ({\R}')
\end{pmatrix}$$ such that $u v = e_{11}$ and $v u = e_{22}.$
\end{remark}

In view of Proposition~\ref{iso}, we immediately have the next:

\begin{cor}\label{isocor} The algebras ${\rm FMat}_{X} (\R)$ and ${\rm FMat}_{X} ({\R}')$ are isomorphic.
\end{cor}

The above corollary is one of the algebraic versions of the
Brown-Green-Rieffel Theorem \cite[Theorem 1.2]{BGR}, and in the
case of countable local units also  follows from \cite[Theorem
2.1]{PereAra}.\\

We shall say that a $G$-graded algebra $\B $ has {\it enough local
units} if for any $g\in G$ the ring ${\D }_g$ has orthogonal local
units. Now it is an easy matter  to derive the next:

\begin{teo}\label{stableteo} Let $G$ be an arbitrary group and
$\B = \bigoplus_{g \in G} \B_g$ be  $G$-graded algebra
 such that
\begin{equation}\label{likeprepr3}  \B_g  \cdot \B_{g\m}  \cdot \B_g
= \B_g  \;\; \forall g \in G.
\end{equation} If $\B $ has enough local units then there exists an infinite  cardinal
$|X|$ and a twisted partial action $\Theta $ of $G$ on $\A ={\rm
FMat}_{X} ({\B}_1)$ such that $${\rm FMat}_{X}(\B ) \cong \A
\stimes_{\Theta} G$$ as graded algebras, the $g$-homogeneous
component of ${\rm FMat}_{X}(\B )$ being ${\rm
FMat}_{X}({\B}_{g}),$ $g\in G.$
\end{teo}
\p By (\ref{likeprepr3}), each ${\B}_g$ is a unital left module
over the ring with orthogonal idempotents ${\D}_g,$ and thus  by
Lemma~\ref{s-unital}, it is torsion-free. Then as in the proof of
Theorem~\ref{criteria2} it follows that $\B $ is homogeneously
non-degenerate.

 Let $X$ be an infinite  set
whose cardinality is bigger than or equal to $|E_g|$ for each
$g\in G,$ where $E_g$ is the set of orthogonal local units of
${\D}_g.$ Add zeros to $E_g$ so that one may assume that $|E_g|=
|X|$ for all $g\in G.$ Write ${\B}' = {\rm FMat}_{X}(\B )$ and $
{\D }'_g = {\B}'_g {\B}'_{g\m}$ $(g\in G).$ Since $\B $ is
homogeneously non-degenerate and satisfies (\ref{likeprepr3}), it
is easy to see that  both of these properties are also satisfied
by ${\B}'.$ Moreover,  for any $g\in G,$  $ {\D }'_g = {\rm
FMat}_{X}({\D }_g)$  and $ {\D }'_g$ has orthogonal local units,
as so does $ {\D }_g.$ Thus our result directly follows from
Theorem~\ref{uv-for-Morita} and Theorem~\ref{criteria2} (or
Theorem~\ref{criteria} in view of Remark~\ref{remark:u,v}). \fim

\begin{remark}\label{countable} As it is observed in the beginning of this section, a countable set of
local units in a ring can be orthogonalized. Thus  the above
theorem can be applied for a  $G$-graded algebra $\B $ which
satisfies (\ref{likeprepr3}) and such that each $\D _g$ has at
most countable set of local units.
\end{remark}

\end{section}

\begin{section}{Examples with uncountable local units}\label{sec:exemplos}

For the concept of a  projective module over a non-necessarily unital ring $\R $ we use the next  definition:
 $_{\R} P$ is said to be {\it projective} if ${\rm Hom} (_{\R} P, \, \underline{\;\;} \; )$
is an exact functor. Evidently, a direct summand of a projective module is again projective.

 We need the following remark, the details can be found in
\cite{Abrams}. We keep denoting by ${\R}\Mod$ the category of all
unital and non-degenerate  ${\R}$-modules.

\begin{remark}\label{fact}
If $\R $ is a ring with local units then for an arbitrary cardinal $\gamma $ the rings
$\R $ and ${\rm FMat}_\gamma(\R )$ are Morita equivalent, and since ${\R }^{(\gamma)}$ and ${\rm FMat}_\gamma(\R
)$ are corresponding objects, we have that $\R $ is projective in  ${\R}\Mod$ if and only if ${\rm
FMat}_\gamma({\R } )$ is  projective  in ${\rm FMat}_\gamma(\R )\Mod$.
\end{remark}

We give the next:

\begin{ex}\label{contraex}  There exist rings $\R $ and ${\R}' $ with local units  such that  $\R $
and ${\R}'$ are Morita equivalent rings; however, ${\rm FMat}_\gamma(\R )$ and ${\rm FMat}_\gamma({\R}' )$ are
not isomorphic for any cardinal $\gamma$.
\end{ex}

\p Let $K$ be a field, and  denote  by $\omega$ the first infinite
ordinal, and by $\alpha $ the least ordinal whose cardinality is
next to that of $\omega ,$ i. e. $|\alpha | = {\aleph}_0^{+}.$

For any element $a\in {\rm RFMat}_\alpha(K)$, we denote by $a(i,j)$ its $\left(i,j\right)$-entry.
As usual,
 $e_{ij}$ stands for the matrix in ${\rm RFMat}_\alpha(K)$ such that $e_{ij}(i,j)=1$ and $0$
otherwise. We write $e_i = e_{ii}.$ Moreover, if $X\subseteq \alpha$ is a subset, then  the ``sum'' $\sum_{i\in
X}e_i$ will be denoted  by $e_X$.

Let $\E $ be the subring of ${\rm RFMat}_\alpha(K)$ formed by the matrices each column of which has
at most
countably many  nonzero entries ($\E=\E_{\alpha,\alpha}$ in the  notation of \cite{RS}).

Set $e=\sum_{i\in\N}e_i=e_\omega$,  $\R =e \E e $ and ${\R}' =\E e \E$. It is  known that $\R $ and ${\R}'$ are
Morita equivalent rings \cite{GS}. Moreover,
\begin{equation}\label{context}
(e\E e, \E e\E, e\E, \E e)
\end{equation} form a Morita context
with surjective bimodule maps given by the multiplication. In
particular, the Morita correspondent of $e \E e$ is $\E e$ and
thus it is a (projective)  generator of $\E e\E $-mod.

It is easily observed that $\R' $ consists of all $\alpha \times \alpha $-matrices with at most countably many
non-zero entries. Clearly $\R$ has local units as it is  a ring with identity, and we check that $\R'$ is also a
ring with local units. Let $N$ be the set of all subsets $X\subset \alpha$ with $\left|X\right|\leq \aleph_0.$
It is easily seen that   $\E e \E=\sum_{X\in N}\E e_X$. Now take $a\in \E e \E$. Then there exist $X, Y \in N$
such that $a=a e_X =e_Y a .$ Setting $Z=X\cup Y$ we have that $a=e_Z a e_Z$; so that $\left\{e_X\right\}_{X\in
N}$, is a set of local units for $\E e \E$.

    If  $\E e \E$ was  a finitely
generated object
in $\E e \E\Mod$ (with respect to the generator $\E e$), then  $\E e \E=\sum \E e \E x_i$ for some
finitely many $x_i \in \E e \E$ which would imply $\E e \E\neq \E e_X$ for some $X\in N.$ Since this is not the case,
$\E e \E$ is not finitely generated $\E e \E$-module.

We shall show that $\E e \E$ is not a projective object in $\E e \E\Mod,$ hence, by Remark
\ref{fact} we will be
done, as $e \E e$ is a unital ring and consequently is projective in $ e \E e \Mod .$

  Observe first that $\Hom_{e\E e}\left(e\E, e\E e\right)\cong
  {\rm RFMat}_\alpha (K)\cdot e$ (notice that it
  is not
$\E e$ as one may suspect). Indeed, let $\f : e\E \to e\E e$ be an
arbitrary homomorphism of left $e\E e$-modules. Evidently $e_{1 t}
\in \E$ and $ e_{1 t} = e e_{1t}$ for any   $t \in \alpha .$ Set
$y_t = \f (e_{1t})$ ($t\in \alpha $) and $y=\sum_{t\in \alpha}
e_{t1}y_t \in {\rm RFMat}_{\alpha }(K)e.$ For any $x \in e\E$ and
$i\in \omega $ write  $e_{i}x = \sum_{j\in F_i} e_i x e_j,$ where
$F_i \subset \alpha $ is a finite subset. Then we have
\begin{align*} & e_i\f (x) = e_{i1} \f (e_{1i} x) = \sum_{j\in F_i}e_{i1}\f (e_{1i}xe_{j1}e_{1j})=
\sum_{j\in
F_i}e_{i1}e_{1i}xe_{j1}\f(e_{1j})=\\ & e_i x \sum_{j\in F_i}e_{j1}y_{j} =  e_ixy.
\end{align*} Since this holds for any $i\in \omega ,$ it follows that $\f (x) = xy$ for any
$x \in e\E .$ Because any
matrix from ${\rm RFMat}_\alpha (K)\cdot e$ obviously determines a homomorphism  $ e\E \to e\E e,$
the claimed
isomorphism follows.

Suppose by contradiction  that $\E e \E$ is a projective object in $\E e \E\Mod$. Then its Morita
correspondent
 $e\E e \E =e \E$ is a projective (and infinitely generated) object in $e\E e\Mod$, and then $e\E$
 has a (infinite) dual
basis $\left\{f_i,x_i\right\}$, where $f_i\in \Hom_{e\E e}\left(e\E, e\E e\right)\cong
{\rm RFMat}_\alpha
(K)\cdot e.$ Hence $f_i=\cdot y_i$, the right multiplication by
$y_i\in {\rm RFMat}_{\alpha}(K)\cdot e$.

We proceed by  showing that there is an infinite sequence
$y_{k_1},\dots, y_{k_n},\ldots \in \{ y_i \}$ and an
element $a\in e\E$, such that $a y_{k_n}\neq 0$ for al $n\in \N$.  For take first  $e_1.$
Clearly $e_1\in e\E$.
There is at least one $y_i$, such that $e_1 y_i\neq 0$. Choose one of them, and denote it by
$y_{k_1}$. Suppose
we have chosen $y_{k_1},\dots, y_{k_n}$, and $e_{11},\dots e_{n i_n}$, such that
\begin{enumerate}
    \item $e_{r i_r} y_{k_r}\neq 0$.
    \item $e_{r i_r} y_{k_s}= 0$, for all $1 \leq r < s \leq n$.
\end{enumerate} By the definition of a dual basis, the set
$\mathcal{U}=\left\{y_i : e_{r i_r}y_i=0,\,1\leq r\leq n \right\}$
can not be empty (in fact, it is  infinite). Pick $y_{k_{n+1}}\in
\mathcal{U},$ and choose $i_{n+1}$ such that $e_{n+1
i_{n+1}}y_{k_{n+1}}\neq 0$. Clearly, all elements in
$\left\{y_{k_n}\right\}$ are distinct.

Now set $a=\sum_{n\in \N}e_{n i_n}.$ Since $e_n a=e_{n i_n}$ we
have that $ay_{k_n}\neq 0$ for all $n\in \N$, which contradicts
the existence of a dual basis.

This means that $e\E$ is not projective and hence $\E e \E$ can
not be projective. \fim \\

\begin{ex}\label{gradedcontraex} Given a group $G,$ which has an element $g$
with $g^2 \neq 1,$ there is a  homogeneously non-degenerate
$G$-graded algebra $\B $ which satisfies (\ref{likeprepr3})  such
that ${\rm FMat}_{\gamma}(\B )$ (with  the $h$-homogeneous
component of ${\rm FMat}_{\gamma}(\B )$ being defined to be ${\rm
FMat}_{\gamma}({\B}_{h}),$ $h\in G$) is not isomorphic as a graded
algebra to a crossed product by a twisted partial action of $G$ on
${\rm FMat}_{\gamma}({\B}_1)$ for any cardinal $\gamma .$
\end{ex}

\p Let $\B$ be the linking algebra of the Morita context
(\ref{context}).  Define a $G$-grading on $\B$ by setting
$$ {\B}_1 = \begin{pmatrix}
e\E e   & 0\\ 0 & \E e\E
\end{pmatrix},\;\; {\B}_g= \begin{pmatrix} 0 & e\E \\
0 &0 \end{pmatrix},\; \; {\B }_{g\m}= \begin{pmatrix} 0  & 0\\
\E e &0\end{pmatrix}, \; \; {\B }_{h}= \begin{pmatrix} 0  & 0\\
0  &0\end{pmatrix}$$ for all $h \in G$ with  $h\neq 1, g, g\m . $
We  have that

$$ {\D }_g = \begin{pmatrix}
e\E e   & 0\\ 0 & 0
\end{pmatrix}, \; \; \; \; {\D }_{g\m } = \begin{pmatrix}
 0  & 0\\ 0 & \E e \E
\end{pmatrix}$$ and  ${\rm FMat}_{\gamma}({\B}_{h}) {\rm FMat}_{\gamma}({\B}_{h\m})
= {\rm FMat}_{\gamma}({\D}_{h}),$ for any $h \in G.$ Moreover, $\B
$ is homogeneously non-degenerate, satisfies (\ref{likeprepr3})
and, consequently both of these properties are verified  by ${\rm
FMat}_{\gamma}({\B}).$  It is directly checked that ${\rm
FMat}_{\gamma}({\D}_g) \cong  {\rm FMat}_{\gamma}( e\E e )$ and
${\rm FMat}_{\gamma}({\D}_{g\m}) \cong {\rm FMat}_{\gamma}( \E e\E
).$ Hence by the previous example ${\rm FMat}_{\gamma}({\D}_g)$ is
not isomorphic to ${\rm FMat}_{\gamma}({\D}_{g\m})$ and thus ${\rm
FMat}_{\gamma}(\B )$ does not satisfy (ii) of
Theorem~\ref{criteria}. \fim

\end{section}

\end{document}